\documentclass[11pt]{amsart}
\usepackage[T1]{fontenc}
\usepackage[utf8]{inputenc}
\usepackage{lmodern}
\usepackage[margin=1in]{geometry}
\usepackage{amsmath,amssymb,mathtools}
\usepackage{booktabs,array}
\usepackage{microtype}
\usepackage[colorlinks=true,linkcolor=blue,citecolor=blue,urlcolor=blue]{hyperref}
\hypersetup{pdftitle={Limits of Generalized Magnitudes of Dissimilarity Spaces},pdfauthor={Hayato Saigo}}
\newtheorem{theorem}{Theorem}[section]
\newtheorem{lemma}[theorem]{Lemma}
\newtheorem{proposition}[theorem]{Proposition}
\newtheorem{corollary}[theorem]{Corollary}
\theoremstyle{definition}
\newtheorem{definition}[theorem]{Definition}
\newtheorem{example}[theorem]{Example}
\theoremstyle{remark}
\newtheorem{remark}[theorem]{Remark}
\newcommand{\one}{\mathbf{1}}
\newcommand{\Mag}{\operatorname{Mag}}
\newcommand{\chiinf}{\chi_\infty}
\newcommand{\mzero}{m_0}
\newcommand{\R}{\mathbb{R}}
\newcommand{\Q}{\mathbb{Q}}
\newcommand{\Z}{\mathbb{Z}}
\DeclareMathOperator{\im}{im}
\DeclareMathOperator{\rank}{rank}
\DeclareMathOperator{\diag}{diag}
\DeclareMathOperator{\adj}{adj}
\DeclareMathOperator{\sgn}{sgn}

\numberwithin{equation}{section}
\title[Limits of generalized magnitudes]{Limits of Generalized Magnitudes of Dissimilarity Spaces}
\author{Hayato Saigo}
\address{ZEN University}
\email{hayato\_saigo@zen.ac.jp}
\date{November 3, 2026}
\begin{document}
\begin{abstract}
For a finite symmetric nonnegative dissimilarity matrix $D=(d_{ij})$ with zero diagonal, we study the high-resolution limit ($t\to \infty$) of $\one^T(e^{-td_{ij}})^+\one$. We prove that this limit always exists in the extended real line and is rational whenever finite, without weightability or invertibility assumptions. At fixed cardinality, a finite rational polyhedral fan governs eventual rank and leading asymptotics. Consequently, only finitely many finite limits occur, and every limit admits an integer-data realization at the same cardinality. In contrast, allowing cardinality to vary, every rational pair of zero-pattern magnitude and high-resolution limit is realized using dissimilarities in $\{0,1,2\}$; adding the level $3$ permits either infinite limit with any rational zero-pattern magnitude. We also classify all limits on at most four points. Moreover, we derive an exact Schur decomposition valid under persistent singularity. Finally, it is shown that positive semidefiniteness at some positive scale is equivalent to equality of dissimilarity rows for every zero-distance pair, in which case the limit counts the zero-distance classes.
\end{abstract}
\maketitle
\section{Introduction}\label{sec:introduction}
The magnitude of a finite similarity matrix extends the counting of points by taking their similarities into account. When a weighting and a coweighting exist, magnitude is their common total; for an invertible matrix $Z$, it is $\one^TZ^{-1}\one$, where $\one$ denotes the all-ones vector. The Moore--Penrose extension
\[
\Mag(Z)=\one^TZ^+\one
\]
is defined for any complex matrix and agrees with the classical value whenever a weighting and a coweighting exist. Chen and Vigneaux explicitly introduced this arbitrary-matrix definition \cite[Definition~4.9]{chen2023}; Akkaya and \"Unl\"u independently developed the extension for finite categories \cite{akkaya2025}. The categorical origins and metric theory of magnitude are developed in \cite{leinster2008,leinster2013}, and its connections with similarity-sensitive diversity appear in \cite{solow1994,leinster2012,meckes2015}.

For a finite metric space, increasing the resolution makes the similarity matrix tend to the identity, and its magnitude tends to the number of points. What remains of this counting interpretation when distinct points can have zero dissimilarity and the triangle inequality is absent? The zero relation need not be transitive. Even when the similarity matrices converge, their pseudoinverses need not do so: their behavior can depend on how the singular limit is approached \cite{stewart1969}. The problem is therefore to distinguish the contribution of the zero pattern from that of the positive dissimilarity levels.

Section~\ref{sec:foundations} fixes the setting. We consider symmetric nonnegative matrices $D=(d_{ij})$ with zero diagonal and the family
\[
A_t(D)=(e^{-td_{ij}}),\qquad
A_\infty(D)=(\mathbf{1}_{\{d_{ij}=0\}}).
\]
The zero-pattern baseline is $\mzero(D)=\one^TA_\infty(D)^+\one$. We call $D$ weightable at infinity when $\one\in\im A_\infty(D)$, but do not impose this condition. Our object is the high-resolution limit
\[
\chiinf(D)=\lim_{t\to\infty}\Mag_D(t),\qquad
\Mag_D(t)=\Mag(A_t(D)).
\]
The baseline is rational, and a nonsingular zero pattern gives $\chiinf(D)=\mzero(D)$. Singular zero patterns require a separate argument.

Section~\ref{sec:asymptotics} supplies that argument without assuming weightability or eventual invertibility. Theorem~\ref{thm:universal-limit} proves that $\chiinf(D)$ always exists in the extended real line and is rational whenever finite. Rank stabilization and Decell's formula express the eventual profile as a quotient of exponential polynomials with integer coefficients. After equal exponents have been combined, their first nonzero terms determine the limit. The same rank analysis gives analyticity away from finitely many exceptional scales. These conclusions precede, and do not depend on, the Schur analysis below.

Section~\ref{sec:cardinality} asks which laws of cardinality survive. Relabelling and positive rescaling preserve the limit, and the additive Cartesian tensor gives an exact product law at every rank. Replication and gluing at a finite positive separation require greater care: their familiar laws hold under the stated weightability assumptions, while explicit counterexamples show why those assumptions matter. Weighted quotient descent and compression by automorphism orbits give further exact reductions.

Section~\ref{sec:fan} makes the arithmetic constraint uniform at fixed cardinality. Theorem~\ref{thm:uniform-fan} constructs a finite rational polyhedral fan on the space of $n$-point data. On each relative interior $C$, the generic rank is fixed and the profile is either eventually zero or has the form
\[
\Mag_D(t)=c_Ce^{-t\ell_C(D)}\bigl(1+O_K(e^{-\delta_Kt})\bigr),
\qquad D\in K\Subset C,
\]
where $c_C\in\Q\setminus\{0\}$ and $\ell_C$ is an integer linear form. The estimate is uniform on compact subsets of $C$. It follows that only finitely many finite limits occur at fixed $n$. Every limit realized by real data also has an integer-data realization at the same cardinality, with the same zero pattern and generic rank; explicit bounds on finite numerators and denominators are obtained as well.

Section~\ref{sec:schur} identifies the source of a correction to the baseline. Splitting the limiting matrix into its range and kernel separates a term tending to $\mzero(D)$ from a kernel contribution. Theorem~\ref{thm:projected-schur} extends this exact decomposition to persistently singular families by first projecting $\one$ onto the current range of $A_t$. The limiting range projector is rational and is constant on each stratum of the same fan. Thus the arithmetic description of the limit is accompanied by a description of how the singular directions contribute to it.

Section~\ref{sec:spectra} applies these results to small cardinalities. It classifies the limits for all nonnegative real data on at most four points and determines sharp thresholds for nonintegral limits, divergence, and strictly positive finite corrections. Exact zero-pattern enumerations through six points and enumerations on five points with dissimilarities restricted to $\{0,1,\ldots,K\}$ for selected positive integers $K\leq6$ are also given.

Section~\ref{sec:realization} allows the cardinality to vary. The attainable joint range is then
\begin{equation}\label{eq:intro-joint-range}
\bigl\{(\mzero(D),\chiinf(D)):D\text{ finite}\bigr\}
=\Q\times\bigl(\Q\cup\{-\infty,+\infty\}\bigr).
\end{equation}
Theorem~\ref{thm:realization} realizes every finite pair using only $\{0,1,2\}$, and Corollary~\ref{cor:joint-range} realizes either infinite limit with any rational baseline using $\{0,1,2,3\}$. The realizations can be weightable at infinity and eventually nonsingular. The construction combines a five-point correction example, binary-pattern substitution, and the gluing law with the matrix-realization theorem of Huang, Tam, and Wu \cite{huang2013}; an elementary binary construction is also provided. Two distinct positive levels are necessary for any finite pair with unequal coordinates, since a finite limit for one-positive-level data must equal its baseline.

Section~\ref{sec:positivity} determines when a class-count interpretation is recovered. Theorem~\ref{thm:psd} proves
\[
\exists t_0>0:\ A_{t_0}\succeq0
\quad\Longleftrightarrow\quad
d_{ij}=0\ \Longrightarrow\ d_{ik}=d_{jk}\quad\text{for all }i,j,k.
\]
These conditions also characterize positive semidefiniteness for every sufficiently large scale. They make the zero relation an equivalence relation compatible with all dissimilarities, and both $\mzero$ and $\chiinf$ count its classes. They do not imply positive semidefiniteness at every positive scale. Outside this class, every positive-scale similarity matrix is indefinite.

Section~\ref{sec:interpretation} returns to the meaning of structural zeros. Strictly positive off-diagonal data form an open dense set on which the limit is $n$, so noninteger limits need not survive small positive perturbations. A joint resolution window recovers the baseline under a rank condition, whereas compactness within a stratum recovers the actual finite limit, including a nonzero correction. These distinctions clarify what the invariant records when observed zeros are retained as structural data or treated as censored positive values, as in the motivating psychophysical setting \cite{kusano2026}.

Appendix~\ref{app:computation} gives exact computational formulas for commensurable dissimilarities. The substitution $x=e^{-ht}$ reduces the exponential problem to a rational-function calculation, including the rank-deficient case. A resolvent constant-term formula provides a second method for symmetric matrices and underlies the finite enumerations.

\section{Dissimilarity spaces and the zero-pattern baseline}\label{sec:foundations}

\begin{definition}[Finite symmetric dissimilarity space]\label{def:dissimilarity}
A finite symmetric dissimilarity space consists of $X=\{1,\ldots,n\}$, $n\geq1$, and $D=(d_{ij})$ with $d_{ij}=d_{ji}\geq0$ and $d_{ii}=0$. No triangle inequality or separation axiom is assumed. For $t\geq0$, define
\begin{equation}\label{eq:magnitude-profile}
 A_t=A_t(D)=(e^{-td_{ij}})_{i,j=1}^n,
 \qquad \Mag_D(t)=\one^TA_t^+\one.
\end{equation}
Its zero-pattern limit is $A_\infty=A_\infty(D)=(\mathbf1_{\{d_{ij}=0\}})$. Set
\[
 \mzero(D)=\one^TA_\infty^+\one.
\]
Call $D$ \emph{weightable at infinity} when $\one\in\im A_\infty$. The \emph{high-resolution limit of generalized magnitude} is
\[
 \chiinf(D)=\lim_{t\to\infty}\Mag_D(t),
\]
when this limit exists in the extended real line.
\end{definition}

The matrices $A_t$ are symmetric but need not be positive semidefinite. Theorem~\ref{thm:universal-limit} establishes existence of the extended limit. Weightability is invoked only where explicitly required or implied by a stated sufficient condition. The data may take arbitrary nonnegative real values; they need not be rational or commensurable.

\begin{lemma}[Quadratic form from a solution]\label{lem:consistent-quadratic}
If $M=M^T\in\R^{n\times n}$, $b\in\im M$, and $Mw=b$, then $b^TM^+b=b^Tw$.
\end{lemma}
\begin{proof}
$M^+b$ solves $Mx=b$, so $w-M^+b\in\ker M=(\im M)^\perp$. Take its inner product with $b$.
\end{proof}

\begin{proposition}[Rationality of the zero-pattern baseline]\label{prop:rational-baseline}
For every $D$, the matrix $A_\infty^+$ has rational entries and $\mzero(D)\in\Q$. If $D$ is weightable at infinity, every solution of $A_\infty w=\one$ satisfies $\mzero(D)=\sum_iw_i$.
\end{proposition}
\begin{proof}
A rational matrix $M$ of rank $r>0$ has a rational full-rank factorization $M=FG$, with $F\in\Q^{m\times r}$ of full column rank and $G\in\Q^{r\times n}$ of full row rank. Then
\[
 M^+=G^T(GG^T)^{-1}(F^TF)^{-1}F^T
\]
is rational. The rank-zero case has $M^+=0$. Apply this to the integer matrix $A_\infty$. Under weightability, use Lemma~\ref{lem:consistent-quadratic}.
\end{proof}

\begin{proposition}[Spectral representation]\label{prop:spectral}
If $\{(\lambda_k(t),v_k(t))\}_{k=1}^n$ is an orthonormal spectral decomposition of $A_t$, then
\begin{equation}\label{eq:spectral}
 \Mag_D(t)=\sum_{\lambda_k(t)\ne0}
 \frac{|\langle\one,v_k(t)\rangle|^2}{\lambda_k(t)}.
\end{equation}
\end{proposition}
\begin{proof}
Insert the spectral decomposition of $A_t^+$ into \eqref{eq:magnitude-profile}.
\end{proof}
Small eigenvalues may cause finite corrections or divergence, with cancellations possible for indefinite matrices.

\begin{proposition}[Nonsingular zero pattern]\label{prop:nonsingular-pattern}
If $A_\infty$ is nonsingular, then $A_t$ is eventually nonsingular and
\[
 \chiinf(D)=\one^TA_\infty^{-1}\one=\mzero(D)\in\Q.
\]
\end{proposition}
\begin{proof}
Continuity of determinant and inversion gives $A_t^{-1}\to A_\infty^{-1}$. Rationality follows from the integer entries of $A_\infty$.
\end{proof}

\section{Rank, regularity, and asymptotics}\label{sec:asymptotics}
This section does not assume weightability or eventual nonsingularity.

\subsection{Exponential polynomials and generic rank}
\begin{definition}[Exponential polynomial]\label{def:exponential-polynomial}
A finite exponential polynomial is a function $f(t)=\sum_{k=1}^m c_ke^{-\lambda_kt}$ with real coefficients and exponents. Its coefficients are integer if $c_k\in\Z$. Equal exponents are combined before leading terms are compared.
\end{definition}

\begin{lemma}[Leading term and zero set]\label{lem:exponential-leading}
For $f\not\equiv0$, written with $\lambda_1<\cdots<\lambda_m$ and $c_k\ne0$,
\[
 f(t)=c_1e^{-\lambda_1t}(1+o(1))\qquad(t\to\infty).
\]
Moreover, $f$ has only finitely many real zeros.
\end{lemma}
\begin{proof}
Factor out $e^{-\lambda_1t}$ for the first claim. For the second, induct on the number of terms. The one-term case has no zeros. For $m>1$, the derivative of $g=e^{\lambda_1t}f$ is a nonzero exponential polynomial with $m-1$ terms. By induction it has finitely many zeros; Rolle's theorem then precludes infinitely many distinct zeros of $g$.
\end{proof}

\begin{definition}[Generic exponential rank]\label{def:generic-rank}
The generic exponential rank $r$ of $A_t$ is the largest order of a minor that is not identically zero as an exponential polynomial.
\end{definition}

\begin{theorem}[Eventual rank stabilization]\label{thm:rank-stabilization}
There is a finite set $F\subset[0,\infty)$ such that $\rank A_t=r$ for every $t\in[0,\infty)\setminus F$. In particular, the rank is eventually constant.
\end{theorem}
\begin{proof}
All minors of order greater than $r$ vanish identically. Choose a nonzero minor of order $r$; Lemma~\ref{lem:exponential-leading} says it has finitely many zeros. Outside those zeros the rank is exactly $r$.
\end{proof}

\subsection{Pseudoinverse formula and finite-scale regularity}
\begin{theorem}[Decell formula \cite{decell1965}]\label{thm:decell}
Let $M\in\R^{m\times n}$ have rank $r>0$, set $B=M^TM$, and write
\[
 \det(\lambda I_n-B)=\lambda^n+a_1\lambda^{n-1}+\cdots+a_r\lambda^{n-r}.
\]
Then $a_r\ne0$ and
\begin{equation}\label{eq:decell}
 M^+=-a_r^{-1}(B^{r-1}+a_1B^{r-2}+\cdots+a_{r-1}I_n)M^T.
\end{equation}
For $r=1$, the parenthesis means $I_n$.
\end{theorem}
\begin{proof}
On a right singular vector $v_i$ with singular value $\sigma_i>0$, the characteristic identity at $\lambda=\sigma_i^2$ gives
\[
 -\frac{\sigma_i^{2r-2}+a_1\sigma_i^{2r-4}+\cdots+a_{r-1}}{a_r}=\sigma_i^{-2}.
\]
Multiplication by $M^T$ yields the corresponding term $\sigma_i^{-1}v_iu_i^T$ of the pseudoinverse. Both sides vanish on $\ker M^T$. Finally, $a_r=(-1)^r\prod_{i=1}^r\sigma_i^2\ne0$.
\end{proof}

\begin{corollary}[Analyticity and fixed-rank continuity]\label{cor:fixed-rank-continuity}
If $M(s)$ is real analytic with constant rank on an interval $I$, then $M(s)^+$ is real analytic on $I$. If $M_k\to M$ and $\rank M_k=\rank M$ for every $k$, then $M_k^+\to M^+$.
\end{corollary}
\begin{proof}
For positive rank, use \eqref{eq:decell}: its numerator and denominator are analytic polynomial expressions in the entries, and the denominator never vanishes. For rank zero the matrices vanish on $I$ (or in the sequence). The same formula proves the sequential assertion; see also \cite{stewart1969}.
\end{proof}

\begin{proposition}[Only finitely many possible discontinuities in scale]\label{prop:finite-discontinuities}
There is a finite $F\subset[0,\infty)$ such that $t\mapsto\Mag_D(t)$ is real analytic on every open interval in $[0,\infty)\setminus F$. In particular, it has at most finitely many discontinuities on $[0,\infty)$.
\end{proposition}
\begin{proof}
Use the rank-exception set from Theorem~\ref{thm:rank-stabilization}. On every interval avoiding it, $A_t$ has constant rank and analytic entries, so Corollary~\ref{cor:fixed-rank-continuity} makes $A_t^+$ analytic. Its quadratic form with $\one$ is analytic. At points outside $F$ the same formula holds in a neighborhood, so discontinuities can occur only in $F$.
\end{proof}

\begin{remark}\label{rem:regularity-dependence}
Proposition~\ref{prop:finite-discontinuities} depends on Theorem~\ref{thm:rank-stabilization} and Corollary~\ref{cor:fixed-rank-continuity}; it does not follow from the eventual representation in Theorem~\ref{thm:universal-limit} alone. At $t=0$, analyticity means restriction of an analytic function in a neighborhood when rank is locally constant.
\end{remark}

\begin{corollary}[No rank increase, no correction]\label{cor:no-rank-increase}
If $\rank A_t=\rank A_\infty$ eventually, then $\Mag_D(t)\to\mzero(D)$.
\end{corollary}
\begin{proof}
Apply the sequential part of Corollary~\ref{cor:fixed-rank-continuity} to $A_t\to A_\infty$.
\end{proof}

\subsection{Eventual representation and the high-resolution limit}
Let $\mathcal E_{\Z}$ be the ring of finite exponential polynomials with integer coefficients and nonnegative real exponents. All entries of $A_t$, its minors, and polynomial matrix expressions in its entries belong to this ring.

\begin{theorem}[Exponential-rational representation and rational finite limits]\label{thm:universal-limit}
For every $D$ in Definition~\ref{def:dissimilarity}, there exist $P,Q\in\mathcal E_{\Z}$ and $T\geq0$ such that $Q(t)\ne0$ and $\Mag_D(t)=P(t)/Q(t)$ for $t\geq T$. Consequently, $\chiinf(D)\in\R\cup\{-\infty,+\infty\}$ always exists; whenever finite, it belongs to $\Q$.
\end{theorem}
\begin{proof}
By Theorem~\ref{thm:rank-stabilization}, $\rank A_t=r$ for all large $t$; $r\geq1$ because every diagonal entry is $1$. Let $G_t=A_t^TA_t$ and write
\[
 \det(\lambda I-G_t)=\lambda^n+a_1(t)\lambda^{n-1}+\cdots+a_r(t)\lambda^{n-r}
\]
on this tail. Apply \eqref{eq:decell} with $Q=a_r$ and
\begin{equation}\label{eq:universal-numerator}
 P=-\one^T\bigl(G_t^{r-1}+a_1G_t^{r-2}+\cdots+a_{r-1}I\bigr)A_t^T\one.
\end{equation}
Both belong to $\mathcal E_{\Z}$, while $Q=(-1)^r\prod_{i=1}^r\sigma_i(t)^2\ne0$ on the tail. If $P\equiv0$, the profile is eventually zero. Otherwise Lemma~\ref{lem:exponential-leading} gives
\[
 P=pe^{-\alpha t}(1+o(1)),\qquad Q=qe^{-\beta t}(1+o(1))
\]
with nonzero integers $p,q$. These coefficients remain integers after terms with equal real exponents have been combined. The quotient tends to $0$, $p/q$, or an infinity with sign $\sgn(p/q)$ according as $\alpha>\beta$, $\alpha=\beta$, or $\alpha<\beta$.
\end{proof}

\begin{corollary}[Rational correction]\label{cor:rational-correction}
If $\chiinf(D)$ is finite, then $\chiinf(D)-\mzero(D)\in\Q$.
\end{corollary}
\begin{proof}
Use Theorem~\ref{thm:universal-limit} and Proposition~\ref{prop:rational-baseline}.
\end{proof}

\begin{remark}\label{rem:correction-dependence}
Corollary~\ref{cor:rational-correction} is an arithmetic conclusion without weightability, even if $A_t$ remains singular. Theorems~\ref{thm:schur} and~\ref{thm:projected-schur} identify the difference as a limiting kernel contribution, including persistent rank deficiency after projection. Section~\ref{sec:fan} strengthens rationality to a finite-value statement when $n$ is fixed.
\end{remark}

\section{Algebraic properties and the analogy with cardinality}\label{sec:cardinality}
The high-resolution invariant has several laws familiar from finite cardinality. Their classical magnitude antecedents include product and distant-union identities \cite[Propositions 1.4.3--1.4.4 and 2.3.6]{leinster2013}; the Moore--Penrose tensor and direct-sum identities are recorded in \cite[Proposition 3.6 and Lemma 4.12]{chen2023}. Their domains of validity differ: the tensor-product law holds without rank assumptions, whereas finite-separation additivity requires control of weightability. This distinction is essential for the Moore--Penrose extension. Call a family \emph{eventually weightable} if $\one\in\im A_t$ for all sufficiently large $t$. This condition concerns the positive-scale family and is distinct from weightability at infinity, which means $\one\in\im A_\infty$. All statements below concern nonempty spaces unless the empty-space convention is explicitly invoked.

\subsection{Normalization, relabelling, and rescaling}
\begin{proposition}[Basic invariances and normalization]\label{prop:invariances}
For every permutation matrix $\Pi$ and every $c>0$,
\[
 \Mag_{\Pi D\Pi^T}(t)=\Mag_D(t),\qquad \Mag_{cD}(t)=\Mag_D(ct).
\]
Consequently, both $\mzero$ and $\chiinf$ are invariant under relabelling and positive rescaling. A one-point space has profile identically $1$. More generally, the all-zero dissimilarity on any nonempty finite set has profile identically $1$. If $D$ has $n$ points and every off-diagonal entry equals $h>0$, then
\[
 \Mag_D(t)=\frac{n}{1+(n-1)e^{-ht}}\quad(t>0),\qquad \chiinf(D)=n.
\]
\end{proposition}
\begin{proof}
Permutation conjugation commutes with the Moore--Penrose inverse and fixes the all-ones vector; rescaling replaces $t$ by $ct$. The zero pattern is unchanged by positive rescaling. For an all-zero dissimilarity, $A_t=J_n$ and $J_n^+=J_n/n^2$. In the last case, the constant vector with entries $1/(1+(n-1)e^{-ht})$ is a weighting; Lemma~\ref{lem:consistent-quadratic} gives the formula.
\end{proof}

\subsection{Replication, quotient weights, and symmetry reduction}
A count of zero-distance classes ignores the addition of another representative of the same class. Generalized magnitude has this property under weightability, but the Moore--Penrose extension does not have it in complete generality. The distinction is already visible in persistently singular families. Failure of unconditional equivalence invariance for generalized categorical magnitude is documented in \cite[Example~2.21]{akkaya2025}; the example below concerns the high-resolution invariant of symmetric dissimilarity spaces.

\begin{proposition}[Exact replication law]\label{prop:replication}
Let $D$ be a dissimilarity matrix on $q$ points. Replace point $i$ by $m_i\geq1$ identical representatives, with dissimilarities pulled back from $D$, and call the resulting matrix $D^{(m)}$. Put $W=\diag(m_1,\ldots,m_q)$ and $s=W^{1/2}\one$. At every scale,
\begin{equation}\label{eq:replication}
 \Mag_{D^{(m)}}(t)=s^T\bigl(W^{1/2}A_t(D)W^{1/2}\bigr)^+s.
\end{equation}
If $\one\in\im A_t(D)$, this equals $\Mag_D(t)$. Uniform replication, $m_1=\cdots=m_q$, preserves the profile at every scale without a weightability assumption. Consequently, it always preserves $\chiinf$, while arbitrary replication preserves $\chiinf$ whenever $A_t(D)$ is eventually weightable.
\end{proposition}
\begin{proof}
Let $P$ be the incidence matrix of the replication map, so that $P^TP=W$, $\one_{\sum_i m_i}=P\one_q$, and $A_t(D^{(m)})=PA_t(D)P^T$. The matrix $Q=PW^{-1/2}$ has orthonormal columns. Therefore
\[
 (PA_tP^T)^+=Q\bigl(W^{1/2}A_tW^{1/2}\bigr)^+Q^T,
\]
which gives \eqref{eq:replication}. If $A_tx=\one_q$, choose $y$ with $P^Ty=x$. Then $(PA_tP^T)y=P\one_q$, and consistency of these symmetric systems gives
\[
 \Mag_{D^{(m)}}(t)=\one_{\sum_i m_i}^Ty=\one_q^Tx=\Mag_D(t).
\]
For $W=mI$, use $(mA_t)^+=m^{-1}A_t^+$ directly in \eqref{eq:replication}.
\end{proof}

\begin{example}[Persistent nonweightability and failure of replication invariance]\label{ex:replication}
Let
\begin{equation}\label{eq:replication-matrix}
 D_*=
 \begin{pmatrix}
 0&1&1&0&0&0\\
 1&0&2&0&1&2\\
 1&2&0&0&1&2\\
 0&0&0&0&1&2\\
 0&1&1&1&0&1\\
 0&2&2&2&1&0
 \end{pmatrix},\qquad x=e^{-t}.
\end{equation}
For $A(x)=(x^{(D_*)_{ij}})$, direct multiplication gives
\[
 A(x)k(x)=0,\qquad k(x)=(-x-1,1,1,-1,1-x,1)^T,
\]
and the leading $5\times5$ principal minor has determinant $-(x-1)^4(x+1)^2$. Hence $A(x)$ has rank $5$ for $0\leq x<1$. Since $\one^Tk(x)=2(1-x)\ne0$ there, the family is nonweightable at every positive scale and at infinity. Writing $B=A(0)$, exact pseudoinversion gives
\[
 B^+\one=\tfrac19(4,2,2,4,2,2)^T,\qquad
 \chiinf(D_*)=\mzero(D_*)=\tfrac{16}{9}.
\]
The last equality also follows from fixed-rank continuity.

Now duplicate only point $2$, placing its copy last, and denote the seven-point dissimilarity by $D_*'$. Its similarity matrices and their limits again have rank $5$. If $B'$ is the new zero-pattern matrix, then
\[
 (B')^+\one=\tfrac1{121}(46,21,20,57,31,31,21)^T,\qquad
 \chiinf(D_*')=\tfrac{227}{121}\ne\tfrac{16}{9}.
\]
Thus arbitrary duplication of an exactly identical point can change even a finite, positive high-resolution value. The exact profiles are
\begin{align*}
 \Mag_{D_*}(t)&=\frac{2(x^3+x^2+6x+8)}{(x+1)(x^2+3)^2},\\
 \Mag_{D_*'}(t)&=\frac{32x^3+32x^2+159x+227}{(x+1)(4x^2+11)^2}.
\end{align*}
These formulas hold for $t>0$; at $t=0$ both matrices have rank $1$ and magnitude $1$.
\end{example}

\begin{proposition}[Weighted descent along a replication quotient]\label{prop:weighted-descent}
For positive point masses $\mu_i$, define
\[
 s_\mu=(\sqrt{\mu_i})_i,\quad
 H_{D,\mu}(t)=\diag(s_\mu)A_t(D)\diag(s_\mu),\quad
 \Mag_{D,\mu}(t)=s_\mu^TH_{D,\mu}(t)^+s_\mu.
\]
Let $\pi:X\twoheadrightarrow Y$ satisfy $d_X(x,x')=d_Y(\pi x,\pi x')$, and give $Y$ the pushed-forward masses $\nu_y=\sum_{\pi x=y}\mu_x$. Then
\begin{equation}\label{eq:weighted-descent}
 \Mag_{D_X,\mu}(t)=\Mag_{D_Y,\nu}(t)\qquad\text{at every scale}.
\end{equation}
Multiplication of all masses by one positive constant does not change the profile. In particular, assigning each of the $m_i$ copies of point $i$ mass $1/m_i$ makes arbitrary replication preserve the original unweighted profile, at every rank.
\end{proposition}
\begin{proof}
Let $P$ be the incidence matrix of $\pi$, and set $Q=\diag(s_\mu)P\diag(s_\nu)^{-1}$. Then $Q^TQ=I$, $s_\mu=Qs_\nu$, and $H_{D_X,\mu}=QH_{D_Y,\nu}Q^T$. Taking the Moore--Penrose inverse through this isometric inclusion proves \eqref{eq:weighted-descent}. Uniform rescaling multiplies $H$ by the scale factor and $s$ by its square root, so it cancels from the quadratic form.
\end{proof}
This is an exact quotient law with specified masses. Unweighted replication pushes counting measure forward to the multiplicities $m_i$, rather than to unit mass at every quotient point; the distinction accounts for Example~\ref{ex:replication}. The weighted construction therefore restores invariance under subdivision of mass without asserting an invariance that the unweighted invariant lacks.

\begin{proposition}[Compression by automorphism orbits]\label{prop:orbit-compression}
Let a finite permutation group act on the points of $D$ by dissimilarity-preserving permutations. Let $O_1,\ldots,O_q$ be its orbits, $m_a=|O_a|$, $P$ their incidence matrix, $W=\diag(m_a)$, and $Q=PW^{-1/2}$. Then
\begin{equation}\label{eq:orbit-compression}
 \Mag_D(t)=s^TH_t^+s,\qquad s=(\sqrt{m_a})_a,\qquad H_t=Q^TA_tQ.
\end{equation}
Explicitly,
\[
 (H_t)_{ab}=\frac1{\sqrt{m_am_b}}\sum_{i\in O_a}\sum_{j\in O_b}e^{-td_{ij}}.
\]
Equivalently, if $L_{ab}(t)=\sum_{j\in O_b}e^{-td_{ij}}$ for any $i\in O_a$, then $H_t=W^{1/2}L(t)W^{-1/2}$. These formulas require neither invertibility nor weightability.
\end{proposition}
\begin{proof}
The orthogonal projection $QQ^T$ onto orbit-constant vectors is the average of the permutation matrices of the group. It commutes with $A_t$, so this subspace and its orthogonal complement reduce $A_t$. Since $\one=Qs$, only the first block contributes to $\one^TA_t^+\one$, proving \eqref{eq:orbit-compression}. Transitivity on each orbit makes $L_{ab}$ independent of $i$ and gives the displayed identities.
\end{proof}
In general, $H_t$ has entries that are sums of exponentials, and its diagonal entries need not equal $1$. The orbit formula is an exact symmetry reduction of the matrix problem; it does not identify the invariant with the unweighted magnitude of a $q$-point dissimilarity space.

\begin{corollary}[Constant row sums and transitive symmetry]\label{cor:constant-row-sum}
Suppose $A_t(D)\one=r_t\one$ at every positive scale. Then
\[
 \Mag_D(t)=\frac{n}{r_t},\qquad \chiinf(D)=\frac{n}{k},
 \qquad k=\#\{j:d_{ij}=0\},
\]
where $k$ is independent of $i$ and includes $j=i$. In particular, these formulas hold whenever the dissimilarity automorphism group acts transitively on the points.
\end{corollary}
\begin{proof}
The common row sum is positive, so $\one/r_t$ is a weighting and its sum is $n/r_t$. Equivalently, $\one$ lies in an eigenspace of the symmetric matrix $A_t$ with nonzero eigenvalue $r_t$; this argument requires no invertibility. Passing to the limit in any row gives $r_t\to\#\{j:d_{ij}=0\}=k$. Transitivity makes all row sums equal.
\end{proof}

\begin{proposition}[Regular zero patterns at every cardinality]\label{prop:transitive-spectrum}
Let $n\geq1$ and $1\leq k\leq n$. A simple $(k-1)$-regular graph on $n$ vertices exists if and only if $n(k-1)$ is even. For every such pair $(n,k)$ and every $h>0$, there is a vertex-transitive dissimilarity matrix with entries in $\{0,h\}$ and profile
\begin{equation}\label{eq:transitive-profile}
 \Mag_D(t)=\frac{n}{k+(n-k)e^{-ht}},\qquad \chiinf(D)=\frac{n}{k}.
\end{equation}
Consequently, the cardinality-$n$ spectrum contains
\[
 \left\{\frac nk:1\leq k\leq n,\ n(k-1)\text{ even}\right\}.
\]
This set is exactly the high-resolution spectrum of the vertex-transitive dissimilarity spaces on $n$ points.
\end{proposition}
\begin{proof}
Necessity is the degree-sum identity. For sufficiency, put $d=k-1$ and take the vertex set $\Z/n\Z$. If $d$ is even, join every vertex to its translates by $\pm1,\ldots,\pm(d/2)$. If $d$ is odd, the parity hypothesis forces $n$ even; use the translates $\pm1,\ldots,\pm((d-1)/2)$ together with the antipodal translate $n/2$. The assumptions $0\leq d\leq n-1$ ensure that these are distinct nonzero translates. These circulant graphs are simple, $d$-regular, and vertex-transitive; the case $d=0$ uses no edges.

Set $d_{ij}=0$ on the diagonal and on graph edges, and set $d_{ij}=h$ on all remaining pairs. Every row of $A_t$ contains $k$ entries equal to $1$ and $n-k$ entries equal to $e^{-ht}$. Corollary~\ref{cor:constant-row-sum} gives \eqref{eq:transitive-profile}. Conversely, the zero-neighbor graph of a vertex-transitive dissimilarity space is regular. Its degree is $k-1$, so the degree-sum identity and Corollary~\ref{cor:constant-row-sum} give the classification within that class.
\end{proof}
For even $n$, this construction realizes every value $n/k$ with $1\leq k\leq n$; for odd $n$, it realizes those with odd $k$. This classifies the transitive case, not the full spectrum. The integer $k$ counts zero-neighbors, not necessarily the size of an equivalence class: noninteger values $n/k$ remain compatible with transitive symmetry.

\subsection{Multiplication: the additive Cartesian tensor}
\begin{proposition}[An exact product law at every rank]\label{prop:tensor}
For dissimilarities $D$ on $X$ and $E$ on $Y$, define $D\boxplus E$ on $X\times Y$ by
\[
 (D\boxplus E)_{(i,a),(j,b)}=d_{ij}+e_{ab}.
\]
Then, at every scale,
\begin{align*}
 A_t(D\boxplus E)&=A_t(D)\otimes A_t(E),\\
 \Mag_{D\boxplus E}(t)&=\Mag_D(t)\Mag_E(t),\\
 \mzero(D\boxplus E)&=\mzero(D)\mzero(E).
\end{align*}
At any fixed scale, the product is weightable if and only if both factors are weightable; the same assertion holds at infinity. If both high-resolution limits are finite, then
\[
 \chiinf(D\boxplus E)=\chiinf(D)\chiinf(E).
\]
The same limit identity holds whenever the corresponding extended-real product is defined, excluding the indeterminate case $0\cdot(\pm\infty)$.
\end{proposition}
\begin{proof}
The exponential of a sum factors entrywise. Singular-value decompositions prove $(A\otimes C)^+=A^+\otimes C^+$, including singular factors, and $\one_{X\times Y}=\one_X\otimes\one_Y$. These identities give the profile and baseline formulas. For weightability, $\im(A\otimes C)=\im A\otimes\im C$. For nonzero vectors $u,v$, membership of $u\otimes v$ in this subspace is equivalent to $u\in\im A$ and $v\in\im C$: apply the orthogonal projection onto $(\im A)^\perp$ in the first factor, and then the analogous projection in the second. Take $u,v$ to be the all-ones vectors. The limit claims follow from the ordinary product rules on their stated domains.
\end{proof}
The tensor is associative and symmetric up to the canonical relabellings of Cartesian products, with the one-point space as unit. These are monoidal identities; no universal property of a categorical Cartesian product is asserted for the additive dissimilarity. If two nonzero tail profiles have leading forms $c_De^{-\lambda_Dt}(1+o(1))$ and $c_Ee^{-\lambda_Et}(1+o(1))$, their product has leading coefficient $c_Dc_E$ and exponent $\lambda_D+\lambda_E$. This follows from the exact product identity and resolves any particular $0\cdot\infty$ case once the leading orders are known.

\subsection{Addition: gluing at a finite positive separation}
For $h>0$, write $D\sqcup_h E$ for the dissimilarity on $X\amalg Y$ that retains the given within-component dissimilarities and assigns the value $h$ to every cross-component pair. This operation is associative and symmetric, up to relabelling, when the same $h$ is used throughout.

\begin{proposition}[The gluing law for weightable blocks]\label{prop:gluing}
At a scale where both $A=A_t(D)$ and $C=A_t(E)$ are weightable, put $a=\Mag_D(t)$, $b=\Mag_E(t)$, and $y=e^{-ht}$. If $1-y^2ab\ne0$, then $A_t(D\sqcup_h E)$ is weightable and
\begin{equation}\label{eq:gluing-profile}
 \Mag_{D\sqcup_h E}(t)=\frac{a+b-2yab}{1-y^2ab}.
\end{equation}
The blocks need not be invertible. If they are invertible, then also
\begin{equation}\label{eq:gluing-determinant}
 \det A_t(D\sqcup_h E)=\det A\det C(1-y^2ab).
\end{equation}
For every $D,E$, without further assumptions,
\[
 \mzero(D\sqcup_h E)=\mzero(D)+\mzero(E).
\]
Weightability at infinity of the glued space is equivalent to weightability at infinity of both factors.
\end{proposition}
\begin{proof}
Choose $Aw=\one_X$ and $Cz=\one_Y$. Then $\one_X^Tw=a$ and $\one_Y^Tz=b$ by Lemma~\ref{lem:consistent-quadratic}. The matrix of the gluing is
\[
 M=\begin{pmatrix}A&y\one_X\one_Y^T\\y\one_Y\one_X^T&C\end{pmatrix}.
\]
With
\[
 \alpha=\frac{1-yb}{1-y^2ab},\qquad \beta=\frac{1-ya}{1-y^2ab},
\]
we have $M(\alpha w,\beta z)^T=\one_{X\amalg Y}$ because $\alpha+yb\beta=1$ and $ya\alpha+\beta=1$. Summing the entries and applying Lemma~\ref{lem:consistent-quadratic} proves \eqref{eq:gluing-profile}. For invertible blocks, Schur elimination and the matrix determinant lemma give \eqref{eq:gluing-determinant}. At infinity, the zero pattern is $A_\infty(D)\oplus A_\infty(E)$; both remaining assertions follow from direct sums.
\end{proof}

\begin{corollary}[Conditional additivity of the limit]\label{cor:additivity}
Suppose $\one\in\im A_t(D)$ and $\one\in\im A_t(E)$ for all sufficiently large $t$, and both $\chiinf(D),\chiinf(E)$ are finite. Then $D\sqcup_h E$ is eventually weightable and
\[
 \chiinf(D\sqcup_h E)=\chiinf(D)+\chiinf(E)\qquad\text{for every }h>0.
\]
Eventual nonsingularity of the factors is sufficient but unnecessary. If the factors are eventually nonsingular, so is the gluing.
\end{corollary}
\begin{proof}
Both scalar profiles are bounded on a tail, so $y^2ab\to0$ and $yab\to0$. Apply \eqref{eq:gluing-profile}, and, for the final assertion, \eqref{eq:gluing-determinant}.
\end{proof}

\begin{example}[Failure of unconditional finite-separation additivity]\label{ex:nonadditivity}
Let $D_*$ be the six-point space of Example~\ref{ex:replication}, and glue a singleton at separation $h>0$. With $y=e^{-ht}$, its similarity matrix is
\[
 M_t=\begin{pmatrix}A_t(D_*)&y\one\\y\one^T&1\end{pmatrix}.
\]
For $t>0$, $\ker A_t(D_*)$ is spanned by $k_t=k(e^{-t})$, with $\one^Tk_t\ne0$. The matrix $M_t$ is nonsingular: in a homogeneous system, multiplication of the first block equation by $k_t^T$ forces the singleton coordinate to vanish; the remaining equations then force the coefficient of $k_t$ to vanish.

For a weighting $(w,z)$, the first block equation multiplied by $k_t^T$ gives $z=1/y$, while the second gives $\one^Tw=(y-1)/y^2$. Thus
\[
 \Mag_{D_*\sqcup_h\{*\}}(t)=\frac2y-\frac1{y^2}
 =2e^{ht}-e^{2ht}\longrightarrow-\infty.
\]
The sum of the component limits is instead $16/9+1=25/9$. Hence finite component limits alone do not ensure finite-separation additivity; a hypothesis such as eventual weightability in Corollary~\ref{cor:additivity} is essential.
\end{example}

For eventually weightable factors, if one or both limits are infinite, the exact formula, rather than formal extended-real addition, controls the result. A useful sufficient condition for preserving a divergent profile $a_t$, when $b_t$ has a finite limit, is $e^{-2ht}a_t\to0$: the denominator tends to $1$, the numerator divided by $a_t$ tends to $1$, and the glued profile is asymptotic to $a_t$.

\subsection{Analogy with cardinality}
On the class of eventually weightable spaces with finite limit, $\chiinf$ is a rational-valued multiplicative invariant for $\boxplus$ and an additive invariant for every fixed positive-separation operation $\sqcup_h$, normalized to $1$ on a singleton. This class is closed under these operations. The tensor and gluing formulas also give the numerical distributive identity for their limits. At the level of dissimilarity spaces, however, $D\boxplus(E\sqcup_h F)$ generally has cross-component dissimilarities $d_{ij}+h$, so it is not identified with $(D\boxplus E)\sqcup_h(D\boxplus F)$.

There is an exact semiring formulation if infinite cross-dissimilarities are explicitly allowed. Set $e^{-t\infty}=0$ for $t>0$ and use disjoint union with cross-dissimilarity $\infty$. Its similarity matrix is a block direct sum. Up to distance-preserving bijections, these extended spaces form a commutative semiring under disjoint union and $\boxplus$, after adjoining the empty space. At each fixed $t>0$, generalized magnitude is a semiring homomorphism to $\R$:
\[
 \Mag_{D\amalg E}(t)=\Mag_D(t)+\Mag_E(t),\qquad
 \Mag_{D\boxplus E}(t)=\Mag_D(t)\Mag_E(t),
\]
with empty-space value $0$ and singleton value $1$. No weightability or rank hypothesis is needed for this matrix-direct-sum statement. The proof of Theorem~\ref{thm:universal-limit} extends to these matrices because the additional entries are identically zero. Consequently, the extended spaces with finite high-resolution limit form a subsemiring, and $\chiinf$ is a semiring homomorphism from it to $\Q$.

\section{Uniform asymptotics and arithmetic finiteness}\label{sec:fan}
The pointwise limit theorem admits a uniform, finite description at each fixed cardinality. The polyhedral mechanism is elementary comparison of monomial weights, as in initial-form constructions \cite[Sections 3 and 5]{maclagan2012}. The issue specific to the present problem is to retain every weight level, including those exposed by exact cancellation, and to treat changes of rank in the same construction.

Put $m=\binom n2$ and identify the data space $\mathcal D_n$ with $\R_{\geq0}^{m}$ using the off-diagonal entries $d=(d_e)_{e\in E}$, where $E$ is the set of unordered pairs of distinct points. For a polynomial $F=\sum_\alpha c_\alpha Y^\alpha$, write $\|F\|_{\mathrm{coef}}=\sum_\alpha|c_\alpha|$.

\begin{lemma}[Universal polynomial data]\label{lem:universal-polynomials}
There are polynomials $P_r,Q_r\in\Z[Y_e:e\in E]$, $1\leq r\leq n$, depending only on $n$, with the following properties:
\begin{enumerate}
\item $\rank A_t(d)\geq r$ if and only if $Q_r(e^{-td})\ne0$.
\item Whenever $\rank A_t(d)=r$, $\Mag_d(t)=P_r(e^{-td})/Q_r(e^{-td})$.
\item Their degrees and coefficient norms satisfy
\[
 \deg Q_r\leq2r,\quad \deg P_r\leq2r-1,\quad
 \|Q_r\|_{\mathrm{coef}}\leq n^{2r},\quad
 \|P_r\|_{\mathrm{coef}}\leq rn^{2r}.
\]
\end{enumerate}
Here $e^{-td}$ denotes the vector with entries $e^{-td_e}$.
\end{lemma}
\begin{proof}
Let $M(Y)$ have diagonal entries $1$ and off-diagonal entry $Y_e$, put $G(Y)=M(Y)^TM(Y)$, and define $a_j(Y)$ by
\[
 \det(\lambda I-G)=\lambda^n+a_1\lambda^{n-1}+\cdots+a_n,
 \qquad a_0=1.
\]
Take
\[
 Q_r=a_r,\qquad
 P_r=-\one^T\left(\sum_{j=0}^{r-1}a_jG^{r-1-j}\right)M^T\one.
\]
Cauchy--Binet gives
\begin{equation}\label{eq:sum-squared-minors}
 (-1)^rQ_r(Y)=\sum_{\substack{I,J\subseteq\{1,\ldots,n\}\\|I|=|J|=r}}
 \det(M(Y)_{I,J})^2.
\end{equation}
For real substitutions, this vanishes exactly when every $r$-minor vanishes. This proves (1), and Decell's formula proves (2). Every entry of $M$ has degree at most $1$ and coefficient norm $1$; every entry of $G$ has degree at most $2$ and coefficient norm at most $n$. The degree assertions follow immediately. Also,
\[
 \|a_j\|_{\mathrm{coef}}\leq\binom nj j!n^j\leq n^{2j}.
\]
For $k\geq0$, expansion of the matrix product yields
\[
 \|\one^TG^kM^T\one\|_{\mathrm{coef}}\leq n^{2k+2}.
\]
Consequently, every one of the $r$ summands defining $P_r$ has coefficient norm at most $n^{2r}$, proving (3).
\end{proof}

\begin{theorem}[Uniform asymptotics]\label{thm:uniform-fan}
For every $n\geq1$, there is a finite rational polyhedral fan $\mathcal F_n$ with support $\mathcal D_n$ such that the following statements hold for the relative interior $C$ of every cone.
\begin{enumerate}
\item The generic exponential rank is a constant $r_C\geq1$. For every compact $K\subset C$, there is $T_K<\infty$ such that $\rank A_t(d)=r_C$ for every $d\in K$ and $t\geq T_K$.
\item Either $\Mag_d(t)=0$ for all sufficiently large $t$, uniformly for $d$ in compact subsets of $C$, or there are a nonzero $c_C\in\Q$ and an integer linear form $\ell_C$ such that
\begin{equation}\label{eq:uniform-asymptotics}
 \Mag_d(t)=c_Ce^{-t\ell_C(d)}\bigl(1+O_K(e^{-\delta_Kt})\bigr),
 \qquad d\in K,\quad t\to\infty,
\end{equation}
for every compact $K\subset C$ and some $\delta_K>0$.
\item The high-resolution limit of generalized magnitude is constant on $C$. In the second alternative it is
\[
 \chiinf(d)=
 \begin{cases}
 0,&\ell_C(d)>0,\\
 c_C,&\ell_C(d)=0,\\
 \sgn(c_C)\,\infty,&\ell_C(d)<0.
 \end{cases}
\]
Each comparison holds either throughout $C$ or nowhere on $C$.
\end{enumerate}
The fan can be chosen using central hyperplanes with integer normal vectors of $\ell^1$-norm at most $4n$.
\end{theorem}
\begin{proof}
Let
\[
 \mathcal A_n=\{\alpha\in\Z_{\geq0}^{m}:|\alpha|_1\leq2n\}.
\]
Use the central hyperplanes $(\alpha-\beta)\cdot d=0$, for distinct $\alpha,\beta\in\mathcal A_n$, together with the coordinate hyperplanes, and restrict their arrangement fan to $\R_{\geq0}^{m}$. The resulting finite rational fan has the asserted support and normal bound. All monomials in Lemma~\ref{lem:universal-polynomials} have exponents in $\mathcal A_n$.

Fix a relative interior $C$. For any one of these polynomials $F$, group its terms according to equality of $\alpha\cdot d$ on $C$. The equality classes and their strict order are independent of $d\in C$. Thus
\begin{equation}\label{eq:grouped-exponents}
 F(e^{-td})=\sum_{j=1}^{s}b_je^{-t\lambda_j(d)},
 \qquad \lambda_1(d)<\cdots<\lambda_s(d)\quad(d\in C),
\end{equation}
where every $\lambda_j$ is the restriction of an integer linear form and every $b_j$ is a fixed integer, namely the sum of the coefficients in its class. Delete all classes with $b_j=0$, not merely the first such class. Either no terms remain, so the specialization is identically zero for every $d\in C$, or its first remaining coefficient is a fixed nonzero integer. Moreover, on every compact $K\subset C$, every remaining gap $\lambda_j(d)-\lambda_1(d)$ with $j>1$ has a positive minimum. Therefore a nonzero specialization satisfies
\begin{equation}\label{eq:uniform-polynomial-leading}
 F(e^{-td})=b_1e^{-t\lambda_1(d)}\bigl(1+O_K(e^{-\delta_Kt})\bigr)
\end{equation}
uniformly on $K$. If only one term remains, the remainder is zero and any $\delta_K>0$ may be used.

Let $r$ be the largest index for which $Q_r(e^{-td})$ is not identically zero on $C$. This index exists because $Q_1$ is the negative sum of the squares of the entries of $M$, whose diagonal entries are $1$. Lemma~\ref{lem:universal-polynomials}(1) and \eqref{eq:uniform-polynomial-leading} show that the rank is at most $r$ at every scale and equals $r$ on a tail uniform over $K$. They also show that this $r$ is the generic exponential rank.

On that tail the profile equals $P_r/Q_r$. If the specialization of $P_r$ is identically zero, this proves the first alternative. Otherwise, let $p,q$ be the first surviving coefficients of $P_r,Q_r$, and let $\alpha,\beta\in\mathcal A_n$ represent the corresponding exponent classes. Divide their uniform expansions, taking the tail sufficiently far that the denominator's relative error is at most $1/2$. This proves \eqref{eq:uniform-asymptotics} with $c_C=p/q$ and $\ell_C(d)=(\alpha-\beta)\cdot d$. The defining arrangement fixes the sign of this difference on $C$, proving (3). When $n=1$, the data space is a point and the claims reduce to $\Mag(t)=1$.
\end{proof}

\begin{corollary}[Finite arithmetic spectrum and integral realization]\label{cor:finite-spectrum}
Fix $n\geq1$. Every finite high-resolution limit of generalized magnitude can be written in lowest terms as $a/b$, where
\[
 a\in\Z,\qquad |a|\leq n^{2n+1},\qquad b\in\Z,\qquad1\leq b\leq n^{2n}.
\]
In particular, the set of finite values at cardinality $n$ is finite. Every high-resolution value realized by nonnegative real dissimilarities, including either infinite value if it occurs, is realized by a nonnegative integer-valued dissimilarity matrix on the same number of points. The integer representative can be chosen with the same zero pattern and generic exponential rank.
\end{corollary}
\begin{proof}
In the nonzero finite case, the limit is $p/q$ for grouped coefficients of $P_r,Q_r$. Grouping cannot increase coefficient norm, and Lemma~\ref{lem:universal-polynomials} gives $|p|\leq rn^{2r}$ and $1\leq|q|\leq n^{2r}$. Reduction of the fraction can only decrease numerator and denominator. The stated bounds follow from $r\leq n$; the zero value is $0/1$.

The relative interior of every nonempty rational polyhedral cone contains a rational point. Clearing denominators gives an integer point in the same relative interior, since it is a cone; for the zero cone take its origin. Apply Theorem~\ref{thm:uniform-fan}. The coordinate hyperplanes in the arrangement ensure equality of zero patterns as well.
\end{proof}

\section{Range--kernel analysis and persistent singularity}
\label{sec:schur}

We first give an inverse-based formula, then extend the decomposition to every
rank by projecting onto the current range. Neither general block formula
assumes weightability. Ordinary Schur elimination is classical; successive
Schur complements also underlie the asymptotic spectral analysis
in~\cite{usevich2026}. Here the objective is the scalar invariant and its
arithmetic limit. Let $R=\im A_\infty$, $N=\ker A_\infty$, and choose matrices
$U,K$ whose columns are respective bases; these need not be orthonormal. Set
$T=[U\ K]$, $u=U^T\one$, $v=K^T\one$, and
\begin{equation}\label{eq:schur-blocks}
\begin{aligned}
B_t&=U^TA_tU,& C_t&=U^TA_tK,& E_t&=K^TA_tK,\\
S_t&=E_t-C_t^TB_t^{-1}C_t,& z_t&=C_t^TB_t^{-1}u,& q_t&=v-z_t.
\end{aligned}
\end{equation}
The last three expressions are defined for sufficiently large $t$, because
$B_\infty=U^TA_\infty U$ is invertible. Norms below are Euclidean vector norms
and their induced operator norms in these fixed coordinates. When
$\ker A_\infty=0$, the kernel block and its quadratic form are absent, and the
formula reduces to ordinary inversion.

\begin{lemma}[The stable range block]\label{lem:stable-range}
$B_\infty$ is invertible, and for every $D$, including those not weightable at
infinity,
\[
\mzero(D)=u^TB_\infty^{-1}u.
\]
\end{lemma}
\begin{proof}
The restriction of $A_\infty$ to its range is bijective, so its bilinear form
in the basis $U$ is nonsingular. Put $w=UB_\infty^{-1}u$. Then $A_\infty w$
is the orthogonal projection of $\one$ onto $R$: it belongs to $R$ and
$U^T(A_\infty w-\one)=0$. The unique solution of this projected equation in
$R$ is $A_\infty^+\one$. Hence
$\one^Tw=u^TB_\infty^{-1}u=\mzero(D)$.
\end{proof}

\begin{theorem}[Exact Schur decomposition]\label{thm:schur}
Assume only that $A_t$ is nonsingular for all sufficiently large $t$. Then
$S_t$ is eventually nonsingular and
\begin{equation}\label{eq:schur-scalar}
\Mag_D(t)=
\underbrace{u^TB_t^{-1}u}_{m_{\mathrm{range}}(t)}+
\underbrace{(v-z_t)^TS_t^{-1}(v-z_t)}_{m_{\mathrm{kernel}}(t)}.
\end{equation}
The range term tends to $\mzero(D)$. Thus $\chiinf(D)$ is finite if and only
if the kernel term has a finite limit $L$, and then
$\chiinf(D)=\mzero(D)+L$.
\end{theorem}
\begin{proof}
The congruence
\[
H_t=T^TA_tT=\begin{pmatrix}B_t&C_t\\C_t^T&E_t\end{pmatrix}
\]
is eventually nonsingular. Its exact block factorization is
\begin{equation}\label{eq:schur-factorization}
H_t=
\begin{pmatrix}I&0\\C_t^TB_t^{-1}&I\end{pmatrix}
\begin{pmatrix}B_t&0\\0&S_t\end{pmatrix}
\begin{pmatrix}I&B_t^{-1}C_t\\0&I\end{pmatrix}.
\end{equation}
Consequently $S_t$ is nonsingular and $A_t^{-1}=TH_t^{-1}T^T$. Use
$T^T\one=(u^T,v^T)^T$ and complete the quadratic form to
obtain~\eqref{eq:schur-scalar}. This argument does not require orthonormal
bases. Finally $B_t^{-1}\to B_\infty^{-1}$, so
Lemma~\ref{lem:stable-range} applies.
\end{proof}

\begin{proposition}[Weightability and the kernel coordinate]
\label{prop:kernel-weightability}
$D$ is weightable at infinity if and only if $v=0$. In that case the kernel
term in~\eqref{eq:schur-scalar} is $z_t^TS_t^{-1}z_t$.
\end{proposition}
\begin{proof}
Since $\im A_\infty=(\ker A_\infty)^\perp$ and the columns of $K$ span the
kernel, $\one\in\im A_\infty$ is equivalent to $K^T\one=v=0$.
\end{proof}

\begin{corollary}[A norm sufficient condition]\label{cor:schur-norm}
Under the hypotheses of Theorem~\ref{thm:schur}, the condition
\[
\|S_t^{-1}\|\,\|q_t\|^2\longrightarrow0
\]
implies $\Mag_D(t)\to\mzero(D)$. Weightability at infinity need not be assumed
separately: the norm condition itself implies it. If $D$ is weightable at
infinity and $C_t=0$ eventually, the same convergence conclusion holds.
\end{corollary}
\begin{proof}
Use $|q_t^TS_t^{-1}q_t|\leq\|S_t^{-1}\|\|q_t\|^2$. If $N\ne0$, then
$C_t\to0$, $E_t\to0$, $B_t^{-1}\to B_\infty^{-1}$, and hence $S_t\to0$,
$z_t\to0$, and $q_t\to v$. Since
\[
\|S_t^{-1}\|\geq\frac1{\|S_t\|}\longrightarrow\infty,
\]
the norm condition forces $v=0$. Apply
Proposition~\ref{prop:kernel-weightability}. If $N=0$, weightability at
infinity is automatic. In the final assertion, $v=z_t=0$.
\end{proof}

The norm estimate discards sign cancellations. Under the eventual
nonsingularity assumption of Theorem~\ref{thm:schur}, a finite limit for a
space that is not weightable at infinity requires cancellation in the
indefinite kernel quadratic form, and cannot be detected by this sufficient
condition. This observation does not extend without qualification to
persistently singular families: fixed-rank continuity can itself give a
finite limit, as in Example~\ref{ex:replication}.

\begin{proposition}[A nonweightable leading obstruction]
\label{prop:nonweightable-obstruction}
Assume the hypotheses of Theorem~\ref{thm:schur} and
$S_t=\sigma_t(S+o(1))$, where $\sigma_t>0$, $\sigma_t\to0$, and $S$ is
nonsingular. Then
\begin{equation}\label{eq:nonweightable-obstruction}
m_{\mathrm{kernel}}(t)
=\sigma_t^{-1}\bigl(v^TS^{-1}v+o(1)\bigr).
\end{equation}
If $v^TS^{-1}v\ne0$, the profile diverges with its sign. If this coefficient
vanishes, the formula alone does not settle the limit.
\end{proposition}
\begin{proof}
$C_t\to U^TA_\infty K=0$ and $B_t^{-1}\to B_\infty^{-1}$, so $z_t\to0$ and
$q_t=v+o(1)$. Also $S_t^{-1}=\sigma_t^{-1}(S^{-1}+o(1))$.
\end{proof}

\subsection{Single-scale lifting under weightability}

\begin{theorem}[Single-scale order-balance asymptotics]
\label{thm:order-balance}
Assume the hypotheses of Theorem~\ref{thm:schur} and weightability at
infinity. Suppose $\rho_t,\sigma_t>0$ tend to zero and
$C_t=\rho_t(C+o(1))$, $S_t=\sigma_t(S+o(1))$ in operator norm, with $S$
invertible. Set $z=C^TB_\infty^{-1}u$. Then
\begin{equation}\label{eq:order-balance}
m_{\mathrm{kernel}}(t)
=\frac{\rho_t^2}{\sigma_t}\bigl(z^TS^{-1}z+o(1)\bigr).
\end{equation}
For $\rho_t=e^{-\alpha t}$ and $\sigma_t=e^{-\beta t}$, put
$c=z^TS^{-1}z$. The correction tends to zero if $2\alpha>\beta$, tends to
$c$ if $2\alpha=\beta$, and diverges with sign $c$ if $2\alpha<\beta$ and
$c\ne0$.
\end{theorem}
\begin{proof}
$B_t^{-1}\to B_\infty^{-1}$ yields $z_t=\rho_t(z+o(1))$, while
$S_t^{-1}=\sigma_t^{-1}(S^{-1}+o(1))$. Weightability at infinity gives
$q_t=-z_t$. Substitute into~\eqref{eq:schur-scalar}.
\end{proof}

\begin{remark}[Multiscale kernel lifting]\label{rem:multiscale}
If $c=0$ in the divergent-order case, later terms must be examined unless
the kernel contribution vanishes identically. Invertibility of $S$ in
Theorem~\ref{thm:order-balance} excludes genuinely multiscale lifting, such
as $S_t\sim\diag(e^{-\beta_1t},e^{-\beta_2t})$ with $\beta_1\ne\beta_2$.
A singular leading coefficient calls for a further splitting of the kernel
and a higher-order analysis. The universal scalar
Theorem~\ref{thm:universal-limit} still determines the eventual scalar
leading order without resolving each kernel direction.
\end{remark}

If some $d_{ij}>0$, let $\delta=\min\{d_{ij}:d_{ij}>0\}$ and
$F_{ij}=\mathbf1_{\{d_{ij}=\delta\}}$. Then
$A_t=A_\infty+e^{-\delta t}F+o(e^{-\delta t})$.

\begin{theorem}[Nondegenerate first kernel lift]\label{thm:first-lift}
Suppose $D$ has a positive entry. Assume weightability at infinity and
invertibility of the compression $P_NF|_N$, where $P_N$ is the orthogonal
projection onto $N$. Then $A_t$ is eventually invertible and
$\Mag_D(t)\to\mzero(D)$; more precisely,
$m_{\mathrm{kernel}}(t)=O(e^{-\delta t})$.
\end{theorem}
\begin{proof}
Put $\varepsilon=e^{-\delta t}$. Then $B_t=B_\infty+O(\varepsilon)$,
$C_t=\varepsilon U^TFK+o(\varepsilon)$, and
$E_t=\varepsilon K^TFK+o(\varepsilon)$. The compression is invertible
exactly when $K^TFK$ is. Hence $S_t=\varepsilon(K^TFK+o(1))$ is invertible,
which makes $A_t$ invertible. Since weightability at infinity gives
$q_t=-z_t=O(\varepsilon)$, the correction is $O(\varepsilon)$.
\end{proof}

\begin{remark}\label{rem:singular-first-lift}
If the first compression is singular, higher dissimilarity levels and
couplings may determine the lift. When $\ker A_\infty=0$, continuity of
ordinary inversion already gives the conclusion.
\end{remark}

\subsection{Persistent rank deficiency and the projected Schur formula}

The range--kernel separation extends to every dissimilarity family once the
forcing vector is projected onto the actual range of $A_t$. This projection
is essential when $\one\notin\im A_t$. The extension uses ordinary
elimination of a consistent linear system, rather than a congruence rule
for the Moore--Penrose inverse.

\begin{proposition}[The limiting range projector]\label{prop:range-projector}
Let $\Pi_t=A_tA_t^+$ be the orthogonal projector onto $\im A_t$, and let $r$
be the eventual rank of $A_t$. Then
\[
\Pi_*:=\lim_{t\to\infty}\Pi_t
\]
exists and is a rational orthogonal projector of rank $r$. Moreover,
\[
\Pi_*A_\infty=A_\infty,\qquad \im A_\infty\subseteq\im\Pi_*.
\]
At fixed cardinality, $\Pi_*$ is constant on the relative interiors of the
cones of the fan $\mathcal F_n$ in Theorem~\ref{thm:uniform-fan}.
\end{proposition}
\begin{proof}
For sufficiently large $t$, Decell's formula expresses every entry of
$A_t^+$, and hence every entry of $\Pi_t$, as a quotient of finite
exponential sums with integer coefficients. Comparing the least
nonvanishing exponent levels shows that each quotient has an extended
limit, rational whenever finite. Entries of an orthogonal projector have
absolute value at most $1$, so the limits are finite. Taking limits in
$\Pi_t^T=\Pi_t$, $\Pi_t^2=\Pi_t$, and $\operatorname{tr}\Pi_t=r$ proves
that $\Pi_*$ is an orthogonal projector of rank $r$. Taking limits in
$\Pi_tA_t=A_t$ gives $\Pi_*A_\infty=A_\infty$.

For the last assertion, use the formal matrix $M(Y)$,
$G(Y)=M(Y)^TM(Y)$, and characteristic coefficients $a_j(Y)$ from the
universal polynomial construction. For each possible rank $r$, set
\[
R_r(Y)=-M(Y)\left(\sum_{j=0}^{r-1}a_j(Y)G(Y)^{r-1-j}\right)M(Y)^T,
\qquad a_0(Y)=1.
\]
On the rank-$r$ tail, $\Pi_t=R_r(e^{-td})/a_r(e^{-td})$. Each entry of
$R_r$ has degree at most $2r\leq2n$, so its support is already contained
in the exponent set used to construct $\mathcal F_n$. On each relative
interior, the exponent coincidences and their order, grouped coefficients,
and generic rank are fixed. The leading coefficient ratio, or the strict
decay that gives a zero limit, is therefore fixed for every projector
entry. Thus the same fan suffices.
\end{proof}

\begin{theorem}[Projected Schur decomposition]\label{thm:projected-schur}
Retain $U,K,B_t,C_t,S_t$ from~\eqref{eq:schur-blocks}, and put
\[
g_t=\Pi_t\one,\qquad \widetilde u_t=U^Tg_t,\qquad
\widetilde v_t=K^Tg_t,\qquad
\widetilde q_t=\widetilde v_t-C_t^TB_t^{-1}\widetilde u_t.
\]
For all sufficiently large $t$, $\widetilde q_t\in\im S_t$ and
\begin{equation}\label{eq:projected-schur}
\Mag_D(t)=\widetilde u_t^TB_t^{-1}\widetilde u_t
+\widetilde q_t^TS_t^+\widetilde q_t.
\end{equation}
No weightability or eventual invertibility assumption is required. The
first term tends to $\mzero(D)$. Consequently the second term has an
extended limit, and
\[
\chiinf(D)=\mzero(D)
+\lim_{t\to\infty}\widetilde q_t^TS_t^+\widetilde q_t.
\]
If $A_t$ is eventually invertible, this reduces to Theorem~\ref{thm:schur}.
\end{theorem}
\begin{proof}
Fix a sufficiently large $t$ so that $B_t$ is invertible, and suppress the
subscript. Since $g\in\im A$, the equation $Ax=g$ has a solution. Writing
$x=Ty$ with $T=[U\ K]$ gives
\[
\begin{pmatrix}B&C\\C^T&E\end{pmatrix}
\begin{pmatrix}y_1\\y_2\end{pmatrix}
=\begin{pmatrix}\widetilde u\\\widetilde v\end{pmatrix}.
\]
Elimination yields $Sy_2=\widetilde q$, proving consistency of the reduced
system. We may therefore choose
\[
y_2=S^+\widetilde q,\qquad
y_1=B^{-1}(\widetilde u-Cy_2).
\]
For every solution $x$ of $Ax=g$, the scalar $g^Tx$ equals $g^TA^+g$:
any two solutions differ by a vector in $\ker A=(\im A)^\perp$. Also
$\one^TA^+\one=g^TA^+g$, because $A^+(I-\Pi)=0$. Hence
\begin{align*}
\Mag_D(t)&=g^TTy\\
&=\widetilde u^TB^{-1}\widetilde u
+(\widetilde v^T-\widetilde u^TB^{-1}C)y_2\\
&=\widetilde u^TB^{-1}\widetilde u+\widetilde q^TS^+\widetilde q.
\end{align*}
This proves the identity without requiring orthonormal bases.

Proposition~\ref{prop:range-projector} implies $\Pi_*U=U$, since the
columns of $U$ lie in $\im A_\infty$. By symmetry,
\[
\widetilde u_t\longrightarrow U^T\Pi_*\one=U^T\one=u.
\]
Thus the first term tends to $u^TB_\infty^{-1}u=\mzero(D)$ by
Lemma~\ref{lem:stable-range}. Theorem~\ref{thm:universal-limit} then gives
the assertion about the second term, including infinite limits. In the
eventually invertible case, $\Pi_t=I$ and $S_t^+=S_t^{-1}$.
\end{proof}

\section{Small-cardinality spectra}\label{sec:small-spectra}\label{sec:spectra}

For a dissimilarity matrix $D$, let $G_0(D)$ be the simple graph whose edges
are the unordered pairs $\{i,j\}$ with $i\ne j$ and $d_{ij}=0$. Thus
$A_\infty=I+A(G_0(D))$. The small-scale limit $t\to0$ has a different
arithmetic behavior; see~\cite{roff2025}. Here all limits are taken as
$t\to\infty$. The following result classifies the limits over all
nonnegative real dissimilarities, with no bound on their positive values.

\begin{theorem}[Complete small-cardinality spectrum]
\label{thm:small-spectrum}
There is no divergent high-resolution limit of generalized magnitude on at most
four points. The exact attainable limit sets at each cardinality are
\begin{center}
\begin{tabular}{cc}
\toprule
$n$ & $\{\chiinf(D):|D|=n\}$\\
\midrule
$1$ & $\{1\}$\\
$2$ & $\{1,2\}$\\
$3$ & $\{1,2,3\}$\\
$4$ & $\{1,4/3,2,3,4\}$\\
\bottomrule
\end{tabular}
\end{center}
More precisely, the joint ranges of $(\mzero,\chiinf)$ are
\begin{align*}
\mathcal J_1&=\{(1,1)\},\\
\mathcal J_2&=\{(1,1),(2,2)\},\\
\mathcal J_3&=\{(1,1),(2,1),(2,2),(3,3)\},\\
\mathcal J_4&=\{(1,1),(4/3,4/3),(2,1),(2,2),\\
&\hspace{3em}(3,1),(3,2),(3,3),(4,4)\}.
\end{align*}
\end{theorem}

We prove the theorem by classifying zero patterns and resolving every
singular pattern. First, the case of one zero pair admits a formula in
arbitrary cardinality.

\begin{proposition}[A single zero pair]\label{prop:single-pair}
Suppose $n\geq2$, $d_{12}=0$, and every other off-diagonal dissimilarity is
strictly positive. If $d_{1j}=d_{2j}$ for every $j\geq3$, then
$\chiinf(D)=\mzero(D)=n-1$. Otherwise put
\[
\rho=\min_{\substack{j\geq3\\d_{1j}\ne d_{2j}}}
\min(d_{1j},d_{2j}),\qquad
\sigma_j=\lim_{t\to\infty}e^{\rho t}
\bigl(e^{-td_{1j}}-e^{-td_{2j}}\bigr).
\]
Each $\sigma_j$ belongs to $\{-1,0,1\}$. With
$k=\sum_{j=3}^n\sigma_j^2\geq1$ and $s=\sum_{j=3}^n\sigma_j$,
\begin{equation}\label{eq:single-pair-limit}
\mzero(D)=n-1,\qquad \chiinf(D)=n-1-\frac{s^2}{k}.
\end{equation}
The matrices $A_t$ are eventually nonsingular in this second case. Every
combination $1\leq k\leq n-2$, $|s|\leq k$, $s\equiv k\pmod2$ is
realizable using only the dissimilarities $0,1,2$.
\end{proposition}
\begin{proof}
Here $A_\infty=J_2\oplus I_{n-2}$, with baseline $n-1$. Use the range basis
$U=[e_1+e_2,e_3,\ldots,e_n]$ and kernel basis $K=e_1-e_2$. Put
$R_t=U^TA_tU$, $C_t=U^TA_tK$, and $u=U^T\one$. Then
\[
R_t\longrightarrow\diag(4,1,\ldots,1),\qquad
C_t=(0,p_3(t),\ldots,p_n(t))^T,
\]
where
\[
p_j(t)=e^{-td_{1j}}-e^{-td_{2j}},\qquad u=(2,1,\ldots,1)^T.
\]
Since $K^TA_tK=0$ and $K^T\one=0$, the scalar Schur complement and its
residual vector are
\[
S_t=-C_t^TR_t^{-1}C_t,\qquad q_t=-C_t^TR_t^{-1}u.
\]
$R_t$ is positive definite for all sufficiently large $t$. If every $p_j$
vanishes, the corresponding two rows are identical and the rank is
eventually $n-1$, so fixed-rank continuity gives the first assertion.
Otherwise $C_t\ne0$ for every positive $t$, and $S_t<0$ eventually,
proving eventual nonsingularity. Moreover,
\[
e^{\rho t}C_t\longrightarrow(0,\sigma_3,\ldots,\sigma_n)^T,
\qquad e^{2\rho t}S_t\longrightarrow-k,
\qquad e^{\rho t}q_t\longrightarrow-s.
\]
Theorem~\ref{thm:schur} gives $\chiinf(D)=(n-1)-s^2/k$, including $s=0$.
To realize the stated values, choose $k$ active coordinates and set
$(d_{1j},d_{2j})=(1,2)$ or $(2,1)$ to prescribe their signs; on inactive
coordinates set both entries equal to $2$. All remaining positive entries
may be set to $1$.
\end{proof}

\begin{lemma}[A clique and one external point]\label{lem:clique-point}
Suppose $G_0(D)=K_r\sqcup K_1$ with $r\geq2$. Write
$u_i(t)=e^{-td_{i,r+1}}$ for the positive cross-dissimilarities. If all
$d_{i,r+1}$ are equal to $h$, then
$\Mag_D(t)=2/(1+e^{-ht})$ and $\chiinf(D)=2$. Otherwise $\Mag_D(t)=1$ for
every $t>0$. In either case $\mzero(D)=2$.
\end{lemma}
\begin{proof}
The matrix has block form $\begin{pmatrix}J_r&u\\u^T&1\end{pmatrix}$.
Writing $s$ for the sum of the first $r$ weights and $z$ for the last
weight, the top equations are $s+u_i z=1$. If $u$ is constant, these reduce,
together with the last equation, to the two-point weighting equations and
give the asserted total. If $u$ is not constant, choose first-block
weights $w$ with $\one^Tw=1$ and $u^Tw=1$, and take $z=0$. Such a choice
exists because $u$ and $\one$ are linearly independent. This exact
weighting has total $1$, so the Moore--Penrose magnitude is $1$. Distinct
positive dissimilarities give distinct similarities at every $t>0$.
\end{proof}

\begin{lemma}[Two zero pairs]\label{lem:two-pairs}
Suppose $G_0(D)=K_2\sqcup K_2$. Then $\mzero(D)=2$ and
$\chiinf(D)\in\{1,2\}$. More explicitly, write
\[
A_t=\begin{pmatrix}J_2&C_t\\C_t^T&J_2\end{pmatrix},\qquad
C_t=\begin{pmatrix}u&v\\w&z\end{pmatrix},\qquad
\delta=u-v-w+z,\quad h=uz-vw,
\]
where the four entries of $C_t$ are exponentials of strictly positive
dissimilarities. If $\delta\not\equiv0$, the limit is $2$. If
$\delta\equiv0$, the limit is $1$ unless all four entries are equal, in
which case it is $2$.
\end{lemma}
\begin{proof}
Direct determinant expansion gives
\[
\det A_t=(h-\delta)(h+\delta),\qquad
\one^T\adj(A_t)\one=2\delta(h-\delta).
\]
We first prove $h=o(\delta)$ whenever $\delta\not\equiv0$. Let $\rho>0$
be the smallest of the four cross-dissimilarities. If the coefficient of
$e^{-\rho t}$ in $\delta$ is nonzero, then $\delta$ has that order while
$h=O(e^{-2\rho t})$. Otherwise the entries at this smallest exponent have
cancelling signs. There are exactly two such entries, one of each sign:
four would make all entries equal and $\delta\equiv0$. The cancelling
entries share a row or a column of $C_t$ and are exactly equal. Direct
substitution gives $h=f\delta$, where $f$ is their common value. Since
$f\to0$, this also proves the claim. Thus $A_t$ is eventually nonsingular
and
\[
\Mag_D(t)=\frac{2\delta}{h+\delta}\longrightarrow2.
\]
If $\delta\equiv0$, linear independence of distinct real exponentials
implies equality of the two multisets of exponents on its positive and
negative sides. Hence either $u=v$ and $w=z$, or $u=w$ and $v=z$. By
interchanging the blocks we may assume $u=v=f$ and $w=z=g$. If $f=g$,
the block totals solve the two-point equations and the profile is
$2/(1+f)$. If $f\ne g$, take the sum of the last two weights to be zero
and choose the first two weights $p,q$ with $p+q=1$, $fp+gq=1$. This is
an exact weighting of total $1$. The dichotomy is proved without assuming
invariance of pseudoinversion under arbitrary duplication.
\end{proof}

\begin{proof}[Proof of Theorem~\ref{thm:small-spectrum}]
If a vertex of $G_0(D)$ is adjacent to every other vertex, its column in
$A_t$ is $\one$ at every scale. The associated unit vector is an exact
weighting, so $\Mag_D(t)=1$. If $A_\infty$ is nonsingular, ordinary
continuity gives $\chiinf(D)=\one^TA_\infty^{-1}\one$, independently of
the positive levels.

For three points the four zero graphs, up to relabeling, are $3K_1$,
$K_2\sqcup K_1$, $P_3$, and $K_3$. The first has limit $3$; the second
has baseline $2$ and limits $1$ or $2$ by Lemma~\ref{lem:clique-point};
the last two have a universal zero vertex and limit $1$. The cases
$n=1,2$ are immediate.

For four points the following table lists all eleven simple graphs up to
isomorphism. The path $P_j$ and cycle $C_j$ have $j$ vertices; the paw is
a triangle with one pendant vertex, and the diamond is $K_4$ with one
edge deleted.
\begin{center}
\begin{tabular}{cccc}
\toprule
$G_0$ & $\rank(I+A(G_0))$ & $\mzero$ & possible $\chiinf$\\
\midrule
$4K_1$ & $4$ & $4$ & $4$\\
$K_2\sqcup2K_1$ & $3$ & $3$ & $1,2,3$\\
$P_3\sqcup K_1$ & $4$ & $2$ & $2$\\
$2K_2$ & $2$ & $2$ & $1,2$\\
$K_3\sqcup K_1$ & $2$ & $2$ & $1,2$\\
$K_{1,3}$ & $4$ & $1$ & $1$\\
$P_4$ & $4$ & $2$ & $2$\\
paw & $3$ & $1$ & $1$\\
$C_4$ & $4$ & $4/3$ & $4/3$\\
diamond & $3$ & $1$ & $1$\\
$K_4$ & $1$ & $1$ & $1$\\
\bottomrule
\end{tabular}
\end{center}
The ranks and baseline values follow by elementary row reduction.
Exhaustiveness follows by edge count: there are respectively
$1,1,2,3,2,1,1$ graph types with $0,1,2,3,4,5,6$ edges. The nonsingular
and universal-vertex observations treat every row except the single-pair,
two-pair, and clique-plus-point rows. These are resolved by
Proposition~\ref{prop:single-pair} and
Lemmas~\ref{lem:clique-point}--\ref{lem:two-pairs}.

All table entries occur. For the single-pair graph,
Proposition~\ref{prop:single-pair} realizes limits $1,2,3$. For
$K_3\sqcup K_1$, choose equal cross-dissimilarities for limit $2$ and
unequal ones for limit $1$. For $2K_2$, take all cross-dissimilarities
equal for limit $2$, and take cross-dissimilarity rows $(1,1)$ and $(2,2)$
for limit $1$. All remaining rows allow arbitrary positive
dissimilarities. The same choices prove the displayed joint ranges.
\end{proof}

\begin{corollary}[Sharp cardinality thresholds]\label{cor:small-thresholds}
Four points are necessary and sufficient for a nonintegral finite limit.
Five points are necessary and sufficient for a divergent limit, and also
for a strictly positive finite correction $\chiinf(D)-\mzero(D)>0$.
\end{corollary}
\begin{proof}
The table and its $C_4$ row establish the first assertion. The table
also excludes divergence and positive corrections on at most four points.
Example~\ref{ex:negative-divergence} supplies a divergent five-point
space. Example~\ref{ex:positive-correction}, with baseline $2$ and limit
$7/3$, supplies a strictly positive finite correction on five points.
Example~\ref{ex:positive-divergence} additionally realizes positive
divergence on five points. None of those examples depends on this
corollary.
\end{proof}

\begin{remark}[A five-point value]\label{rem:five-point-eleven-thirds}
For $n=5$, Proposition~\ref{prop:single-pair} gives the one-zero-pair
limit set $\{1,2,3,11/3,4\}$. For example, take the only zero pair to be
$\{1,2\}$ and set
\[
(d_{13},d_{14},d_{15})=(1,1,2),\qquad
(d_{23},d_{24},d_{25})=(2,2,1),
\]
with all dissimilarities among $3,4,5$ equal to $1$. Then
$(\mzero,\chiinf)=(4,11/3)$. This classifies one zero pattern on five
points, not the full five-point spectrum.
\end{remark}

The zero-pattern baseline can be enumerated completely for small
cardinalities. Write
\[
\mathcal B_n=
\{\one^TB^+\one:B\in\R^{n\times n},\ B=B^T,\ B_{ii}=1,\ B_{ij}\in\{0,1\}\}.
\]
Then we have the following:
\begin{equation}\label{eq:baseline-enumeration}
\begin{aligned}
\mathcal B_1&=\{1\},\qquad\mathcal B_2=\{1,2\},\qquad\mathcal B_3=\{1,2,3\},\\
\mathcal B_4&=\{1,\tfrac43,2,3,4\},\\
\mathcal B_5&=\{1,\tfrac43,\tfrac75,\tfrac53,2,\tfrac73,3,4,5\}.
\end{aligned}
\end{equation}
For six points,
\begin{equation}\label{eq:baseline-six}
\begin{aligned}
\mathcal B_6=\{&0,1,\tfrac65,\tfrac43,\tfrac75,\tfrac{10}7,
\tfrac32,\tfrac53,\tfrac{12}7,\tfrac{16}9,\\
&2,\tfrac73,\tfrac{12}5,\tfrac83,3,\tfrac{10}3,4,5,6\}.
\end{aligned}
\end{equation}
The enumeration counts are shown below. Here a labelled zero pattern is a graph on the fixed vertex set $\{1,\ldots,n\}$, with distinct edge sets counted separately. The column ``Labelled'' counts all such graphs, whereas the other columns count isomorphism classes with the indicated properties.
\begin{center}
\begin{tabular}{rrrrr}
\toprule
$n$ & Patterns & Labelled & Singular & Nonweightable\\
\midrule
1 & 1 & 1 & 0 & 0\\
2 & 2 & 2 & 1 & 0\\
3 & 4 & 8 & 2 & 0\\
4 & 11 & 64 & 6 & 0\\
5 & 34 & 1024 & 19 & 0\\
6 & 156 & 32768 & 88 & 1\\
\bottomrule
\end{tabular}
\end{center}
In particular, every zero pattern on at most five points is weightable.
On six points, the unique nonweightable isomorphism type is $B=I+A(G)$,
where $G$ is the tree consisting of two adjacent vertices, each joined to
two leaves. Its baseline is $16/9$. 
We enumerate all subsets of the set of unordered pairs of distinct vertices. Vertex permutations act on these subsets, and we choose one representative from each orbit. We verify that the resulting orbits are pairwise disjoint and cover all labelled zero patterns.

For a positive integer $K$, we next restrict the dissimilarities to the integers from $0$ to $K$ and enumerate all possible high-resolution limits on five points. This finite computation must be distinguished from a classification allowing arbitrary nonnegative real dissimilarities.

Define
\[
\mathcal S_5^{[K]}=
\{\chiinf(D):D\in\mathcal D_5,\ d_{ij}\in\{0,1,\ldots,K\}\}.
\]
The exact enumeration gives
\begin{equation}\label{eq:alphabet-enumeration}
\begin{aligned}
\mathcal S_5^{[1]}&=\mathcal B_5,\\
\mathcal S_5^{[2]}&=\mathcal S_5^{[1]}\cup\{\tfrac{15}7,\tfrac{11}3\},\\
\mathcal S_5^{[3]}&=\mathcal S_5^{[2]}\cup\{0,\tfrac52,+\infty,-\infty\},\\
\mathcal S_5^{[4]}&=\mathcal S_5^{[3]}\cup\{\tfrac95\},\\
\mathcal S_5^{[5]}&=\mathcal S_5^{[6]}=\mathcal S_5^{[4]}.
\end{aligned}
\end{equation}
The conclusion for $K=5$ also follows immediately by inclusion between
the computed $K=4$ and $K=6$ spectra. The absence of further values through
$K=6$ is not asserted to classify arbitrary integer or real data on five
points. In particular, integer realization at fixed cardinality does not
provide an upper bound of $6$ on the entries of an integer representative.

For each singular zero-pattern representative, we assign to each nonedge of its zero graph an integer between $1$ and $K$, and check every such assignment.
For nonsingular zero patterns, the limit is already the inverse sum. The numbers of singular-pattern
assignments checked are
\begin{center}
\begin{tabular}{lrrrrr}
\toprule
$K$ & 1 & 2 & 3 & 4 & 6\\
\midrule
assignments & 19 & 1391 & 33871 & 376789 & 12485059\\
\bottomrule
\end{tabular}
\end{center}
These are representative-assignment counts, not counts of distinct
labelled matrices: automorphisms may cause repeated assignments, but
cannot omit any labelled matrix up to permutation. The exact
rank-sensitive algorithm, including persistent singularity and
nonweightability, is proved in Appendix~\ref{app:computation}.

\begin{example}[The smallest cardinality permitting a zero limit]
\label{ex:zero-limit}
Let
\[
D=\begin{pmatrix}
0&1&0&1&1\\
1&0&2&1&2\\
0&2&0&3&2\\
1&1&3&0&0\\
1&2&2&0&0
\end{pmatrix}.
\]
For $x=e^{-t}$, direct determinant identities give
\[
\det A_t=x^4(1-x)^4(2-x^2),\qquad
\one^T\adj(A_t)\one=x^6(1-x)^4.
\]
Thus $A_t$ is nonsingular for every $t>0$ and
\[
\Mag_D(t)=\frac{x^2}{2-x^2}>0,\qquad \chiinf(D)=0.
\]
The zero-pattern matrix is permutation-equivalent to
$J_2\oplus J_2\oplus I_1$, so $\mzero(D)=3$. The entire baseline is
cancelled by the limiting kernel contribution. Together with the complete
real-data classification in Theorem~\ref{thm:small-spectrum}, this proves
that five is the smallest cardinality at which a zero limit can occur.
\end{example}

\begin{example}[Negative divergence on five points]
\label{ex:negative-divergence}
Let
\[
D_- =\begin{pmatrix}
0&1&0&1&1\\
1&0&1&1&3\\
0&1&0&3&2\\
1&1&3&0&0\\
1&3&2&0&0
\end{pmatrix}.
\]
Its zero-pattern matrix is again permutation-equivalent to
$J_2\oplus J_2\oplus I_1$, so $\mzero(D_-)=3$ and it is weightable at
infinity. Writing $x=e^{-t}$, exact polynomial identities give
\begin{align*}
\det A_t&=x^6(1-x)^4(x^2+2x+2),\\
\one^T\adj(A_t)\one
&=x^4(1-x)^4(x^2+2x+2)(x^2+2x-2).
\end{align*}
The determinant is nonzero for every $t>0$, and hence
\[
\Mag_{D_-}(t)=1+\frac2x-\frac2{x^2},\qquad \chiinf(D_-)=-\infty.
\]
All dissimilarities lie in $\{0,1,2,3\}$. Together with
Theorem~\ref{thm:small-spectrum}, this also proves that five points are
necessary and sufficient for negative divergence.
\end{example}

\section{Realization of arbitrary rational baselines and limits}\label{sec:realization}

The zero-pattern baseline need not determine the high-resolution limit. We now
make their independence precise. The matrix-realization ingredient is classical:
Huang, Tam, and Wu proved that every rational number is the sum of the entries
of the inverse of the adjacency matrix of a connected nonsingular simple
graph~\cite{huang2013}. We first adapt this result to unit-diagonal zero
patterns. The resulting binary realization is not claimed as a new inverse-sum
theorem. We then give the examples and operations needed for independent
realization of the baseline and the limit.

\begin{lemma}[Binary realization]\label{lem:binary-realization}
For every $c\in\Q$, there is a symmetric nonsingular matrix $B_c$ with entries
in $\{0,1\}$ and all diagonal entries equal to $1$, such that
\[
  \one^T B_c^{-1}\one=c.
\]
Consequently, the binary dissimilarity matrix $D_c=J-B_c$ satisfies
$\mzero(D_c)=\chiinf(D_c)=c$, and $A_t(D_c)$ is eventually nonsingular.
\end{lemma}
\begin{proof}
If $c=1$, take $B_c=[1]$. Otherwise put $a=c/(c-1)$, so $a\ne1$.
By the theorem of Huang, Tam, and Wu, there is a symmetric nonsingular
$\{0,1\}$-matrix $H$ with zero diagonal and $\one^TH^{-1}\one=a$.
Set $B_c=J-H$. The matrix determinant lemma gives
\[
  \det B_c=(-1)^n\det H\,(1-a)\ne0.
\]
If $w=H^{-1}\one$, then $(J-H)w=(a-1)\one$, whence
$\one^TB_c^{-1}\one=a/(a-1)=c$. Finally, $A_t(D_c)\to B_c$,
so continuity of inversion proves the assertion about $D_c$.
\end{proof}

\begin{remark}[An elementary construction]\label{rem:elementary-realization}
There is also a direct construction independent of the general graph-realization
theorem. For a nonsingular symmetric binary matrix with unit diagonal, write
$\mu(B)=\one^TB^{-1}\one$. Such matrices are closed under direct sum and
Kronecker product, with
\[
  \mu(B\oplus C)=\mu(B)+\mu(C),\qquad
  \mu(B\otimes C)=\mu(B)\mu(C).
\]
The matrix
\[
B_-=
\begin{pmatrix}
1&0&0&1&0&0&0\\
0&1&1&1&1&0&0\\
0&1&1&0&0&1&0\\
1&1&0&1&0&1&1\\
0&1&0&0&1&1&0\\
0&0&1&1&1&1&0\\
0&0&0&1&0&0&1
\end{pmatrix}
\]
has determinant $1$, as follows by row reduction, and
\[
 B_-(-4,2,-3,5,-3,2,-4)^T=\one.
\]
Thus $\mu(B_-)=-5$ and $N=B_-\oplus I_4$ has order $11$ and
$\mu(N)=-1$. For an integer $q\ge2$, put $k=2q-1$ and
\[
 C_q=\begin{pmatrix}I_k&J_k\\J_k&I_k\end{pmatrix}.
\]
Its eigenvalues are $1$ with multiplicity $2k-2$, and $1+k,1-k$, so it is
nonsingular. Its row sum is $k+1$, and hence
$\mu(C_q)=2k/(k+1)=2-1/q$. Therefore
\[
 F_q=C_q\oplus N\oplus N,\qquad P_q=N\otimes F_q
\]
realize $-1/q$ and $1/q$, respectively. Repeated direct sums realize every
rational with denominator at least $2$; identity matrices and direct sums of
$N$ realize positive and negative integers, and $[1]\oplus N$ realizes zero.
This gives explicit finite matrices without a claim of size optimality.
\end{remark}

\begin{example}[A positive correction on five points]\label{ex:positive-correction}
Let
\[
E=\begin{pmatrix}
0&1&2&0&0\\
1&0&0&1&2\\
2&0&0&1&1\\
0&1&1&0&1\\
0&2&1&1&0
\end{pmatrix}.
\]
Its zero-pattern matrix has rank $4$ and a weighting
$(1,1,0,0,0)^T$, so $\mzero(E)=2$. With $x=e^{-t}$, direct determinant
expansion gives
\begin{align*}
 \det A_t(E)&=-x^2(1-x)^4(x+1)(x+3),\\
 \one^T\adj(A_t(E))\one&=x^2(1-x)^4(x^2-2x-7).
\end{align*}
Consequently, $A_t(E)$ is nonsingular for every $t>0$, and
\[
 \Mag_E(t)=\frac{7+2x-x^2}{(x+1)(x+3)},\qquad
 \bigl(\mzero(E),\chiinf(E)\bigr)=\left(2,\frac73\right).
\]
All dissimilarities belong to $\{0,1,2\}$. The correction is $1/3>0$.
\end{example}

\begin{proposition}[Binary-pattern substitution]\label{prop:pattern-substitution}
Let $B$ be a symmetric nonsingular binary matrix of order $N$ with unit
diagonal, let $E$ be a dissimilarity matrix on $m$ points, and let $h>0$.
Define $\mathcal P_h(B,E)$ on $\{1,\ldots,N\}\times\{1,\ldots,m\}$ by
\[
 \mathcal P_h(B,E)_{(i,a),(j,b)}=
 \begin{cases}
 E_{ab},&B_{ij}=1,\\
 h,&B_{ij}=0.
 \end{cases}
\]
Put $c=\one_N^TB^{-1}\one_N$, $f_t=\Mag_E(t)$, and $y=e^{-ht}$.
At a scale where $A=A_t(E)$ is invertible and $1-yf_t\ne0$,
\begin{align}
 \det A_t(\mathcal P_h(B,E))
 &= (\det B)^m(\det A)^N(1-yf_t)^{N-1}
       \bigl[1+y(c-1)f_t\bigr].\label{eq:pattern-det}
\end{align}
If also $1+y(c-1)f_t\ne0$, then
\begin{equation}
 \Mag_{\mathcal P_h(B,E)}(t)
 =\frac{cf_t}{1+y(c-1)f_t}.\label{eq:pattern-profile}
\end{equation}
If $A_t(E)$ is eventually nonsingular and $\chiinf(E)$ is finite, the
substitution is eventually nonsingular and
\[
 \mzero(\mathcal P_h(B,E))=c\,\mzero(E),\qquad
 \chiinf(\mathcal P_h(B,E))=c\,\chiinf(E).
\]
If $E$ is weightable at infinity, so is the substitution.
\end{proposition}
\begin{proof}
The entrywise definition gives a symmetric nonnegative dissimilarity matrix
with zero diagonal, and its similarity matrix is
\[
 M=B\otimes(A-yJ_m)+yJ_N\otimes J_m.
\]
Let $H=A-yJ_m$. The matrix determinant lemma and the rank-one inverse formula
give
\[
 \det H=\det A(1-yf_t),\qquad
 \one_m^TH^{-1}\one_m=\frac{f_t}{1-yf_t}.
\]
Thus $M$ is a rank-one update of $B\otimes H$, whose inverse-entry sum is
$cf_t/(1-yf_t)$. A second application of the same formulas gives
\eqref{eq:pattern-det} and \eqref{eq:pattern-profile}.
If $f_t$ has a finite limit, both $1-yf_t$ and $1+y(c-1)f_t$ tend to $1$,
which proves eventual nonsingularity and the limit formula.

Since $h>0$, the zero-pattern matrix is $B\otimes A_\infty(E)$, where
$A_\infty(E)=(\mathbf1_{\{E_{ab}=0\}})_{a,b}$. The Kronecker pseudoinverse
identity gives its baseline $c\,\mzero(E)$. If $A_\infty(E)w=\one_m$, then
$(B^{-1}\one_N)\otimes w$ is a weighting of $B\otimes A_\infty(E)$,
proving weightability at infinity.
\end{proof}

\begin{proposition}[One positive dissimilarity level]\label{prop:one-level}
Suppose every dissimilarity belongs to $\{0,h\}$ for some $h>0$, and let
$B$ be the zero-pattern matrix. Put $x=e^{-ht}$, so
$A_t=(1-x)B+xJ$.
If $D$ is weightable at infinity and $m=\mzero(D)$, then, whenever
$1-x+xm\ne0$,
\[
 \Mag_D(t)=\frac{m}{1-x+xm}.
\]
In particular, $\chiinf(D)=\mzero(D)$. If $D$ is not weightable at
infinity, then $\Mag_D(t)=1/x=e^{ht}$ for every $t>0$, and
$\chiinf(D)=+\infty$.
\end{proposition}
\begin{proof}
In the first case, choose $Bw=\one$, so $\one^Tw=m$.
Then $A_tw=(1-x+xm)\one$, and dividing by this nonzero scalar gives an
exact weighting. The denominator tends to $1$ as $t\to\infty$.
In the second case, symmetry and $\one\notin\im B$ give
$k\in\ker B$ with $\one^Tk\ne0$. Since
$A_tk=x(\one^Tk)\one$, the vector $k/[x(\one^Tk)]$ is an exact
weighting of total $1/x$. In both cases the total of a weighting equals the
Moore--Penrose magnitude.
\end{proof}

\begin{theorem}[Independent realization over two positive levels]\label{thm:realization}
For every pair $(a,b)\in\Q^2$, there is a finite symmetric dissimilarity
matrix $D$ with entries in $\{0,1,2\}$ such that
\[
 \mzero(D)=a,\qquad \chiinf(D)=b.
\]
The matrix $D$ can be chosen weightable at infinity and with $A_t(D)$
invertible for all sufficiently large $t$. Two distinct positive levels are
necessary to realize a finite pair with $a\ne b$: every dissimilarity
matrix taking values in $\{0,h\}$, where $h>0$, whose high-resolution
limit is finite satisfies $\chiinf(D)=\mzero(D)$.
\end{theorem}
\begin{proof}
Let $E$ be the five-point space in Example~\ref{ex:positive-correction}.
For the prescribed pair put
\[
 r=3(b-a),\qquad s=7a-6b.
\]
Choose $B_r$ and $D_s=J-B_s$ as in Lemma~\ref{lem:binary-realization}, and set
\[
 D=\mathcal P_2(B_r,E)\sqcup_1 D_s.
\]
Proposition~\ref{prop:pattern-substitution} and the gluing law
(Proposition~\ref{prop:gluing} and Corollary~\ref{cor:additivity}) show that
$D$ is weightable at infinity and eventually nonsingular, with
\[
 \mzero(D)=2r+s=a,\qquad
 \chiinf(D)=\frac73r+s=b.
\]
The substitution uses only the levels already in $E$, since its cross level
is $2$; the final gluing adds only the level $1$. Thus all entries belong to
$\{0,1,2\}$. If $B_r$ and $B_s$ have orders $N_r,N_s$, the construction
has $5N_r+N_s$ points. When $a=b$, the binary realization $D_a$ alone
suffices. The necessity assertion follows from Proposition~\ref{prop:one-level}:
with only one positive level, failure of weightability at infinity forces
$\chiinf(D)=+\infty$, while weightability gives
$\chiinf(D)=\mzero(D)$.
\end{proof}

\begin{example}[Positive divergence on five points]\label{ex:positive-divergence}
Let
\[
F_+=\begin{pmatrix}
0&0&1&1&1\\
0&0&2&3&1\\
1&2&0&0&2\\
1&3&0&0&1\\
1&1&2&1&0
\end{pmatrix}.
\]
Its zero-pattern matrix is $J_2\oplus J_2\oplus I_1$, so it is weightable
at infinity and $\mzero(F_+)=3$. For $x=e^{-t}$, determinant expansion gives
\begin{align*}
 \det A_t(F_+)&=-x^5(1-x)^4(x+2),\\
 \one^T\adj(A_t(F_+))\one&=x^4(1-x)^4(x^2-4).
\end{align*}
Hence $A_t(F_+)$ is nonsingular for every $t>0$, and
\[
 \Mag_{F_+}(t)=\frac2x-1,\qquad \chiinf(F_+)=+\infty.
\]
Together with Theorem~\ref{thm:small-spectrum}, this also shows that five
points are necessary and sufficient for positive divergence.
\end{example}

\begin{corollary}[Attainability range]\label{cor:joint-range}
As finite dissimilarity spaces vary, their attainable baseline--limit pairs are
exactly
\[
 \bigl\{(\mzero(D),\chiinf(D)):D\text{ finite}\bigr\}
 =\Q\times\bigl(\Q\cup\{-\infty,+\infty\}\bigr).
\]
Every pair can be realized by a matrix with entries in $\{0,1,2,3\}$ that
is weightable at infinity and whose similarity matrix is invertible for all
sufficiently large positive scales.
\end{corollary}
\begin{proof}
Rationality of the zero-pattern baseline and Theorem~\ref{thm:universal-limit}
give inclusion in the displayed set. Theorem~\ref{thm:realization} realizes
every finite pair with the stated properties.

Let $F_+$ be Example~\ref{ex:positive-divergence}, and let $F_-$ be the
negative-divergence space in Example~\ref{ex:negative-divergence}. Both have
entries in $\{0,1,2,3\}$, are weightable at infinity, and have nonsingular
similarity matrices at every positive scale. With $x=e^{-t}$,
\[
 \mzero(F_+)=\mzero(F_-)=3,\qquad
 \Mag_{F_+}(t)=\frac2x-1,\qquad
 \Mag_{F_-}(t)=1+\frac2x-\frac2{x^2}.
\]
For a prescribed $a\in\Q$, take the binary realization $D_{a-3}$ from
Lemma~\ref{lem:binary-realization}, and form
\[
 D_+=F_+\sqcup_1D_{a-3},\qquad
 D_-=F_-\sqcup_2D_{a-3}.
\]
Their zero-pattern matrices are direct sums of weightable blocks, with baseline
$3+(a-3)=a$, and their entries still belong to $\{0,1,2,3\}$.

Put $b_t=\Mag_{D_{a-3}}(t)$, which is bounded and tends to $a-3$.
For $D_+$ the cross-similarity is $y=x$, so
$y^2\Mag_{F_+}(t)b_t=O(x)$. For $D_-$ it is $y=x^2$, so
$y^2\Mag_{F_-}(t)b_t=O(x^2)$. In both cases the extra determinant factor
in Proposition~\ref{prop:gluing} tends to $1$, proving eventual
nonsingularity. The same proposition gives
\[
 \Mag_{D_+}(t)=\frac2x+O(1)\longrightarrow+\infty,\qquad
 \Mag_{D_-}(t)=-\frac2{x^2}+\frac2x+O(1)\longrightarrow-\infty.
\]
This realizes either infinite limit with any prescribed rational baseline.
\end{proof}

\section{Positive semidefiniteness}\label{sec:positivity}

A class-count interpretation of the high-resolution limit of generalized magnitude requires
additional structure; finiteness and nonnegativity of the dissimilarities alone
do not supply it. If $A_t\succeq0$, then $A_t^+\succeq0$ and therefore
$\Mag_D(t)\ge0$. The stronger structural conclusion below characterizes the
existence of even one positive semidefinite positive-scale similarity matrix.

\begin{theorem}[Pseudometric case]\label{thm:pseudometric}
If $d$ is a pseudometric and $C_1,\ldots,C_q$ are its zero-distance classes,
let $\bar d$ be the induced metric and
$\bar A_t=(e^{-t\bar d(C_a,C_b)})_{a,b}$. For all sufficiently large $t$,
\begin{equation}
 \Mag_D(t)=\one_q^T\bar A_t^{-1}\one_q\longrightarrow q.
 \label{eq:pseudometric-limit}
\end{equation}
In particular, a metric on $n$ points has limit $n$.
\end{theorem}
\begin{proof}
The triangle inequality makes $\bar d$ well defined. For the class-incidence
matrix $P$, one has $A_t=P\bar A_tP^T$ and $\one_n=P\one_q$.
Since $\bar A_t\to I_q$, it is eventually invertible. Choose $w$ with
$P^Tw=\bar A_t^{-1}\one_q$; such a vector exists because $P$ has full
column rank. Then $A_tw=\one_n$, and the total of this weighting gives
\[
 \Mag_D(t)=\one_n^Tw=\one_q^T\bar A_t^{-1}\one_q\longrightarrow q.
\]
\end{proof}

\begin{theorem}[Positive-scale rigidity]\label{thm:psd}
For a finite dissimilarity matrix $D$, the following conditions are equivalent:
\begin{enumerate}
\item $A_{t_0}$ is positive semidefinite for some $t_0>0$;
\item $d_{ij}=0$ implies $d_{ik}=d_{jk}$ for every $k$;
\item $A_t$ is positive semidefinite for every sufficiently large $t$.
\end{enumerate}
Under these conditions, zero dissimilarity is an equivalence relation, $d$ is
constant on pairs of its classes, and, writing $q$ for their number,
\[
 \chiinf(D)=\mzero(D)=q.
\]
No triangle inequality is needed. If $q>1$ and $\bar\delta$ is the smallest
positive dissimilarity in the quotient, positive semidefiniteness holds whenever
$t>\log(q-1)/\bar\delta$.
\end{theorem}
\begin{proof}
Assume (1). If $d_{ij}=0$, then
$(e_i-e_j)^TA_{t_0}(e_i-e_j)=0$. Positive semidefiniteness implies
$A_{t_0}(e_i-e_j)=0$. Thus columns $i,j$ agree; symmetry and injectivity
of the real exponential at $t_0>0$ prove (2).

Assume (2). Reflexivity and symmetry of the zero relation hold by definition.
If $d_{ij}=d_{jk}=0$, equality of rows $i,j$ gives $d_{ik}=d_{jk}=0$,
proving transitivity. Equality of rows also makes the quotient dissimilarity
$\bar D$ well defined. With class-incidence matrix $P$,
$A_t=P\bar A_tP^T$ and $\bar A_t\to I_q$. Hence $\bar A_t$ is positive
definite for all large $t$, proving (3). More quantitatively, every
off-diagonal quotient entry is at most $e^{-t\bar\delta}$, and strict diagonal
dominance gives positive definiteness when
$(q-1)e^{-t\bar\delta}<1$. If $q=1$, then $A_t=J$ for all $t$.
The implication (3)$\Rightarrow$(1) is immediate.

The weighting argument in Theorem~\ref{thm:pseudometric} applies to this
compatible quotient and gives $\chiinf=q$. Moreover, the zero-pattern matrix
$A_\infty(D)=PP^T$ is a direct sum of all-ones blocks. Assigning weight $1/|C|$
to each point in a class $C$ gives a weighting with total $q$, proving
$\mzero=q$.
\end{proof}

The qualification ``for sufficiently large $t$'' in Theorem~\ref{thm:psd}
is essential. For example, consider the three-point dissimilarity given by
$d_{12}=d_{23}=1$ and $d_{13}=5$. There are no off-diagonal zeros, so
condition~(2) holds. With $x=e^{-t}$,
\[
 \det A_t=(1-x^5)(1+x^5-2x^2).
\]
Since $1+e^{-5t}-2e^{-2t}=-t+O(t^2)$ as $t\downarrow0$, this determinant
is negative for all sufficiently small positive $t$. Thus $A_t$ is indefinite
at those scales, although it is positive definite for all sufficiently large
$t$. In particular, the equivalent conditions of Theorem~\ref{thm:psd} do
not imply positive semidefiniteness at every positive scale.

\begin{corollary}[Indefiniteness at every positive scale]\label{cor:indefinite}
If some pair with $d_{ij}=0$ has distinct dissimilarity rows, then $A_t$ is
indefinite for every $t>0$. In particular, this holds whenever
$\chiinf(D)\notin\{1,\ldots,n\}$ or $\chiinf(D)\ne\mzero(D)$.
\end{corollary}
\begin{proof}
Theorem~\ref{thm:psd} excludes positive semidefiniteness at every positive
scale. Each $A_t$ also has a positive eigenvalue because
$\operatorname{tr}A_t=n>0$. The numerical conditions contradict the theorem's
conclusion and therefore force incompatible zero rows.
\end{proof}

\section{Interpretation and structural zeros}\label{sec:interpretation}

Let $\mathcal D_n$ be the space of symmetric nonnegative $n\times n$ matrices
with zero diagonal, identified with $[0,\infty)^{n(n-1)/2}$.

\begin{proposition}[The generic separated locus]\label{prop:separated-locus}
The set
\[
 \mathcal D_n^{\mathrm{sep}}
 =\{D\in\mathcal D_n:d_{ij}>0\text{ for }i\ne j\}
\]
is relatively open and dense in $\mathcal D_n$. For every $D$ in this set,
$A_\infty(D)=I_n$ and $\chiinf(D)=n$.
\end{proposition}
\begin{proof}
Strict positivity of finitely many coordinates defines a relatively open set.
Adding arbitrarily small positive entries at off-diagonal zeros proves density.
Since $A_t(D)\to I_n$, continuity of inversion gives the limit.
\end{proof}

\begin{example}[Discontinuity at a structural-zero matrix]\label{ex:discontinuity}
Let $D_0$ have four points, with zero dissimilarity on the edges of the cycle
$C_4$ and dissimilarity $1$ on its two diagonals. Its zero-pattern matrix has
constant row sum $3$, and $\chiinf(D_0)=4/3$ by the constant-row-sum formula.
Replace each zero off-diagonal entry of $D_0$ by $\varepsilon>0$, obtaining
$D_\varepsilon$. Then $D_\varepsilon\to D_0$, whereas
$\chiinf(D_\varepsilon)=4$ for every $\varepsilon>0$. Thus
$D\mapsto\chiinf(D)$ is discontinuous at $D_0$.
\end{example}

\begin{proposition}[Structural-zero resolution window]\label{prop:resolution-window}
Let $D\in\mathcal D_n$ have at least one positive entry, and put
$\delta=\min\{d_{ij}:d_{ij}>0\}$. Suppose $D^{(k)}\in\mathcal D_n$,
$\varepsilon_k\downarrow0$, $t_k\to\infty$, and
\begin{equation}
 \begin{gathered}
 0\le d_{ij}^{(k)}\le\varepsilon_k\quad(d_{ij}=0),\qquad
 d_{ij}^{(k)}\ge\delta/2\quad(d_{ij}>0),\\
 t_k\varepsilon_k\longrightarrow0.
 \end{gathered}\label{eq:resolution-window}
\end{equation}
Then $A_{t_k}(D^{(k)})\to A_\infty(D)$. If $A_\infty(D)$ is nonsingular, or more
generally if $\rank A_{t_k}(D^{(k)})=\rank A_\infty(D)$ for all sufficiently large
$k$, then
\[
 \Mag_{D^{(k)}}(t_k)\longrightarrow\mzero(D).
\]
\end{proposition}
\begin{proof}
At a zero entry, $0\le t_kd_{ij}^{(k)}\le t_k\varepsilon_k\to0$, so the
exponential tends to $1$. At a positive entry it is at most
$e^{-t_k\delta/2}\to0$. Thus the matrices converge. Fixed-rank continuity
of the Moore--Penrose inverse gives the assertion; nonsingularity is the
full-rank special case.
\end{proof}

\begin{remark}[Noncommuting limits]\label{rem:noncommuting-limits}
For each fixed separated $D_\varepsilon$, the strict high-resolution limit
is $n$. A joint regime $t\to\infty$, $\varepsilon\to0$ with
$t\varepsilon\to0$ may instead recover $\mzero(D)$ under the rank
condition in Proposition~\ref{prop:resolution-window}. Matrix convergence
alone gives no conclusion about the magnitude profile when the approximating
matrices have higher rank than their limit. In particular, that proposition
recovers the baseline $\mzero(D)$, not a corrected limit
$\chiinf(D)\ne\mzero(D)$ such as $7/3$ in
Example~\ref{ex:positive-correction}; it does not establish robust recovery
of that corrected asymptotic value.
\end{remark}

\begin{corollary}[Recovery within a stratum]\label{cor:stratum-recovery}
Let $C$ be the relative interior of a cone in the fan of
Theorem~\ref{thm:uniform-fan} on which $\chiinf$ is finite, and let
$K\subset C$ be compact. If $D_k\in K$ and $t_k\to\infty$, then
\[
 \Mag_{D_k}(t_k)\longrightarrow\chiinf(C).
\]
This conclusion includes strata on which $\chiinf\ne\mzero$.
\end{corollary}
\begin{proof}
The eventual-zero alternative is immediate. If $\ell_C=0$ on $C$, the
uniform estimate in Theorem~\ref{thm:uniform-fan} gives convergence to
$c_C$. If $\ell_C>0$ on $C$, compactness gives $\min_K\ell_C>0$, and the
same estimate gives convergence to zero. These exhaust the finite-limit cases.
\end{proof}

For pseudometrics the invariant counts zero-distance classes. On the open dense
separated locus of Proposition~\ref{prop:separated-locus}, the limit is $n$;
fractional structural-zero values are not stable under arbitrary positive
perturbations. The resolution-window result recovers $\mzero(D)$ under its
rank condition, whereas Corollary~\ref{cor:stratum-recovery} recovers the actual
limit under its stratum constraint. For general dissimilarities,
the high-resolution limit of generalized magnitude is a singular-perturbation invariant:
$\mzero(D)$ records the zero pattern, while a finite nonzero Schur correction
records how positive levels lift its kernel. Theorem~\ref{thm:universal-limit}
gives an extended limit, rational whenever finite, even when weightability or
eventual invertibility fails. Statistical use should distinguish structural
zeros from small, noisy positive observations; positive and negative
small-eigenvalue contributions may cancel.

In the motivating psychophysical protocol~\cite{kusano2026}, individual
pairwise ratings were recorded on a finite eight-point response scale and
transformed to a finite dissimilarity scale including zero. Distinct stimuli
can therefore receive an exact observed zero, and pairwise judgments do not
enforce transitivity. Retaining that category as a structural observed zero
leads to the zero-pattern analysis above. Treating it instead as a censored
latent value in $(0,\varepsilon]$ leads to the separated limit $n$ for every
fixed positive completion; Proposition~\ref{prop:resolution-window} describes
a joint finite-resolution regime recovering $\mzero(D)$ when its rank
hypothesis holds. Averaging individual ratings may remove exact zeros, so the
choice of individual or pooled data changes the estimand. The general limit
is accordingly a descriptor whose interpretation requires assumptions about
the data and the zero relation.

\appendix
\section{Commensurable dissimilarities and exact computation}\label{app:computation}

\begin{definition}[Commensurable dissimilarities]\label{def:commensurable}
A dissimilarity matrix $D$ is \emph{commensurable} if $d_{ij}=h k_{ij}$ for some $h>0$ and a nonnegative integer matrix $K=(k_{ij})$. Every rational dissimilarity matrix is commensurable. Put $x=e^{-ht}$ and
\[
 M(x)=(x^{k_{ij}})\in\Z[x]^{n\times n}.
\]
\end{definition}

\begin{proposition}[Polynomial generic rank]\label{prop:polynomial-rank}
If $r=\rank_{\Q(x)}M(x)$, then $\rank M(x)=r$ for every sufficiently small $x>0$.
\end{proposition}
\begin{proof}
Every minor of order greater than $r$ vanishes identically, whereas a nonzero $r$-minor has only finitely many positive roots.
\end{proof}

\begin{theorem}[Bordered determinant formula for the high-resolution tail]\label{thm:bordered-determinant}
If $\det M(x)\not\equiv0$, then for every sufficiently small $x>0$,
\begin{equation}\label{eq:bordered-determinant}
 \Mag_D(t)=\one^T M(x)^{-1}\one
 =-\frac{\det\begin{pmatrix}M(x)&\one\\\one^T&0\end{pmatrix}}{\det M(x)}\in\Q(x).
\end{equation}
\end{theorem}
\begin{proof}
Apply Proposition~\ref{prop:polynomial-rank} and the block determinant identity.
\end{proof}

\begin{theorem}[Rational pseudoinverse in the generic rank-deficient case]\label{thm:rational-pseudoinverse}
Suppose $0<r=\rank_{\Q(x)}M(x)<n$. Put $G(x)=M(x)^TM(x)$ and write
\[
 \det(\lambda I-G(x))=\lambda^n+a_1(x)\lambda^{n-1}+\cdots+a_r(x)\lambda^{n-r}.
\]
On the small positive $x$ tail,
\begin{equation}\label{eq:rational-pseudoinverse}
 M(x)^+=-\frac{G(x)^{r-1}+a_1(x)G(x)^{r-2}+\cdots+a_{r-1}(x)I}{a_r(x)}M(x)^T,
\end{equation}
and hence $\one^T M(x)^+\one\in\Q(x)$.
\end{theorem}
\begin{proof}
Combine Proposition~\ref{prop:polynomial-rank} with Theorem~\ref{thm:decell}. The coefficients and numerator are polynomials, and $a_r$ is nonzero on the relevant tail.
\end{proof}

\begin{proposition}[Exact one-sided limit of a rational function]\label{prop:rational-limit}
For $f(x)=P(x)/Q(x)\in\Q(x)$ with $P,Q$ nonzero, put
\[
 \gamma=\operatorname{ord}_0P-\operatorname{ord}_0Q,
 \qquad c=\frac{\operatorname{lc}_0(P)}{\operatorname{lc}_0(Q)},
\]
where $\operatorname{lc}_0$ denotes the coefficient of the lowest nonzero power. Then $f(x)=cx^\gamma(1+o(1))$ as $x\to0^+$. Its limit is $0$ if $\gamma>0$, $c$ if $\gamma=0$, and $\sgn(c)\infty$ if $\gamma<0$.
\end{proposition}
\begin{proof}
Factor the least powers of $x$ from the numerator and denominator.
\end{proof}

Thus one may choose $h$, form $M(x)$, determine its generic rank, apply Theorem~\ref{thm:bordered-determinant} or~\ref{thm:rational-pseudoinverse}, and use Proposition~\ref{prop:rational-limit}. An identically zero numerator gives limit zero. These one-variable calculations are exact specializations of Theorem~\ref{thm:universal-limit}, not prerequisites for it.

For symmetric matrices, a second exact implementation can avoid matrix powers in~\eqref{eq:rational-pseudoinverse}. It also treats nonweightable and persistently singular families.

\begin{proposition}[Resolvent constant-term formula]\label{prop:resolvent-constant-term}
Let
\[
 d(x,\lambda)=\det(M(x)+\lambda I)=\sum_{j=0}^n d_j(x)\lambda^j,
 \qquad
 q(x,\lambda)=\one^T\adj(M(x)+\lambda I)\one
              =\sum_{j=0}^{n-1}q_j(x)\lambda^j.
\]
Let $k$ be the smallest index for which $d_k(x)\not\equiv0$, and set $q_{-1}=0$ and any coefficients beyond the displayed degrees equal to zero. Then, on the small positive $x$ tail,
\begin{equation}\label{eq:resolvent-constant-term}
 \one^T M(x)^+\one
 =\frac{q_k(x)d_k(x)-q_{k-1}(x)d_{k+1}(x)}{d_k(x)^2}.
\end{equation}
\end{proposition}
\begin{proof}
Fix a small positive $x$ at which the generic rank is attained. Since $M(x)$ is real symmetric, its spectral decomposition gives
\[
 \frac{q(x,\lambda)}{d(x,\lambda)}
 =\one^T(M(x)+\lambda I)^{-1}\one
 =\frac{\|P_{\ker M(x)}\one\|^2}{\lambda}
   +\one^T M(x)^+\one+O(\lambda).
\]
Thus the resolvent quadratic form has at most a simple pole at zero, and its constant term is the required magnitude. Expanding the quotient after factoring $\lambda^k$ from $d$ gives~\eqref{eq:resolvent-constant-term}. The identity holds on an interval of positive $x$ and therefore as a rational-function identity. For $k=0$ it reduces to $q_0/d_0$.
\end{proof}

All determinant coefficients in these procedures can be computed using integer arithmetic. In an enumeration with dissimilarities restricted to $\{0,1,\ldots,K\}$ for a positive integer $K$, nonsingular zero patterns require only their inverse sums; singular zero patterns require testing the positive-level assignments. To enumerate zero patterns, identify their edge sets with subsets of the unordered pairs of distinct vertices and partition these subsets into orbits under vertex permutations. Selecting one representative from every orbit and testing every assignment from $\{1,\ldots,K\}$ to its nonedges covers every allowed dissimilarity matrix up to simultaneous permutation of rows and columns. Repeated assignments related by automorphisms do not affect completeness.

\section*{Acknowledgements}
The author thanks Kyoko Kusano and Naotsugu Tsuchiya for the collaboration reported in part in \cite{kusano2026}. This research was supported by JSPS KAKENHI Grant-in-Aid for Transformative Research Areas (A), ``Qualia Structure: Bridging the Gap between Subjective Experience and Objective Science'' (JP23H04830).
\section*{Declaration of generative AI and AI-assisted technologies in the manuscript preparation process}
In accordance with the recommendation of the Leiden Declaration \cite{leiden2026} on Artificial Intelligence and Mathematics, the author discloses the use of ChatGPT sol 6, Claude Fable, and ChatGPT astra to assist in checking and revising the manuscript. The author takes full responsibility for the final manuscript, including the correctness of its mathematical arguments and results and the accuracy and completeness of its references.
\nocite{*}
\bibliographystyle{amsplain}
\bibliography{references}
\enlargethispage{3\baselineskip}
\end{document}